%% file: Amain.tex
\documentclass[oneside]{amsart}

\usepackage{subfiles}	
\usepackage[utf8]{inputenc}
\usepackage[english]{babel}
\usepackage{amsmath,amsthm,mathtools}
\usepackage{amsfonts}
\usepackage{amssymb}
\usepackage{graphicx}
\usepackage[]{geometry}
\usepackage{tikz}
\usepackage{tikz-3dplot,pgfplots}
\usepackage{color}
\usepackage{enumerate}
\usepackage[all]{xy}
\usepackage[normalem]{ulem}
\usepackage{todonotes}
\usepackage{url}
\usepackage{cleveref}

 \newtheorem{theorem}{Theorem}[section]
\newtheorem*{maintheorem}{Main Theorem}
  \newtheorem{corollary}[theorem]{Corollary}

 \newtheorem{lemma}[theorem]{Lemma}
  \newtheorem{proposition}[theorem]{Proposition}

   \theoremstyle{definition}
\newtheorem{definition}[theorem]{Definition}

\theoremstyle{remark}    
    \newtheorem{remark}{Remark}[theorem]
     \newtheorem{example}[theorem]{Example}

 \usepackage{enumitem}
    
\begin{document}

\title{The Word Problem is Solvable for 3-free Artin groups in quadratic time }
\author[R. Blasco]{Rub\'en Blasco-Garc{\'\i}a}
\address{Rub\'en Blasco-Garc{\'\i}a: Departamento de Matemáticas, Universidad de Zaragoza, Spain}
\email{rubenb@unizar.es}
\author[M. Cumplido]{María Cumplido}
\address{María Cumplido: Instituto de Matemáticas de la Universidad de Sevilla (IMUS) and Departmento de Álgebra,
Facultad de Matemáticas, Universidad de Sevilla, Spain}
\email{cumplido@us.es}
\author[R. Morris-Wright]{Rose Morris-Wright}
\address{Rose Morris-Wright: Middlebury College, Middlebury Vermont, USA}
\email{rmorriswright@middlebury.edu}

\begin{abstract}
We give a {quadratic-time explicit and computable algorithm} to solve the word problem for Artin groups that do not contain any relations of length 3. Furthermore, we prove that, given two geodesic words representing the same element, one can obtain one from the other  by using a set of homogeneous relations that never increase the word length. 
\end{abstract}
\maketitle
\section{Introduction} \label{Section: Introduction}

\subfile{Introduction}

\subsection*{Acknowledgements}
Rub\'en Blasco thanks Sarah Rees for the time he spent in Newcastle under her supervision, that was crucial for the existence of this paper. 
Mar\'ia Cumplido was supported by the research project PID2022-138719NA-I00, financed by MCIN/AEI/10.13039/501100011033/FEDER, UE. She also thanks US-1263032 and P20\_01109 financed by Junta de Andaluc\'ia, the Research Program ``Braids" of ICERM (Providence, RI) and the grant Ram\'on y Cajal 2021. All authors thank Derek Holt and Sarah Rees for reading a first version of this paper and helping to improve it. They are also grateful to the anonymous referee for their careful reading of the manuscript and for the helpful suggestions that improved its clarity.

\section{Background} \label{Section: Background}

\subfile{Background}

\section{Rightward Reducing Sequences} \label{Section: RSS}

\subfile{RRS}

\section{Main Results} \label{Section: MainResults}

\subfile{Results}

\section{Proof of key propositions} \label{Section:Propositions}

\subfile{Props}

\section{Proof of technical lemmas} \label{Section:lemmas}

\subfile{Lemmas}

\bibliographystyle{alpha}
\bibliography{references.bib}

\end{document}

%% file: Introduction.tex
Artin (or Artin--Tits) groups are finitely presented groups whose relations are written in the form $sts\cdots=tst\cdots$, where both words in the equality have the same length. Jacques Tits first defined these groups in the 60's, as a generalization of braid groups, and they first appear in the literature in  \cite{Bourbaki}. Since then, group theorists have focused on classical group theoretic questions in Artin groups--- such as the conjugacy problem or the word problem. Results for these questions are known for particular types of Artin groups. However, these remain open questions in the most general setting. 

\medskip

Hereunder is a summary of articles that solve the word problem for some subfamilies of Artin groups:

\begin{itemize}

\item FC-type Artin groups, including braid groups, spherical-type Artin groups and RAAGs. \cite{Artin1947,Garside1969, Adyan, Deligne1972, Brieskorn1972, Epstein1992, Elrifai1994, Altobelli1998, Altobelli2000}.

\item 2-dimensional Artin groups \cite{AppelShupp,Appel,Chermak};

\item Sufficiently large Artin groups, including large and extra-large. \cite{Holt2012,Holt2013};

\item Euclidean Artin groups \cite{Digne2006,Digne2012, McCammond2017}.

\end{itemize}
\medskip

In this article, we give an explicit algorithm to solve the word problem for a big family of Artin groups: the Artin groups that do not contain any relation of length 3 (also called braid relations), that is, no relations $sts=tst$ for any two generators $s$ and $t$. We will call these groups \emph{3-free} Artin groups. Note that commuting relations of the form $st=ts$ are permitted. {The techniques that we use are a non-trivial generalization of the techniques created by Holt and Rees in their previously mentioned papers dealing with sufficiently large Artin groups.}

\begin{theorem}
There is an algorithm that solves the word problem for 3-free Artin groups {whose complexity is quadratic with respect of length of the word}.
\end{theorem}

The collection of 3-free Artin groups overlaps with the types listed above, but also includes new groups for which the problem was previously unsolved. The intersection of 3-free type Artin groups with FC-type and Euclidean Artin groups is very small: just free and direct products of Artin groups with two generators that are not braid groups (including RAAGs). An Artin group is 2-dimensional if, for every three generators, the three respective relations have lengths $m_1$, $m_2$ and $m_3$ and satisfy $\frac{1}{m_1}+\frac{1}{m_2}+\frac{1}{m_3}\leq 1$. All 3-free Artin groups containing a 3-generator subgroup such that $m_1=2,\, m_2=2$ are not 2-dimensional. Sufficiently large Artin groups are those such that for any triple of generators, if $m_1=2$, $m_2\neq 2,\infty$, then $m_3=\infty$. Again, the 3-free Artin groups with some triple such that $m_1=2$, $m_2=2$, $m_3\neq \infty$ are not sufficiently large. Notice also that 3-free Artin groups include the class of even Artin groups, which are the ones that only contain relations of even length. 

The family of 3-free Artin groups is a subset of the set of locally reducible Artin groups. In \cite{Charney2004}, Charney shows that locally reducible Artin groups satisfy the K($\pi$,1) conjecture, implying that these groups are torsion-free. This also implies that they have trivial center as long as they do not factor as direct products (See \cite{Jankiewicz2022}). In addition, Jankeiwicz and Schreve (\cite{jankiewicz2022rightangled}) have shown that locally reducible Artin groups satisfy a generalized version of the Tits conjecture.  However before our work, there was very little known about the algorithmic properties of 3-free Artin groups.  For example, many of the types of Artin groups listed above  are known to be biautomatic (including finite type \cite{Charney1992}, RAAGS \cite{VanWyk1994,Hermiller1995}, and extra-large type \cite{Peifer1996}). This paper aims to provide tools for proving other properties of 3-free Artin groups. 

 Our algorithm for solving the word problem (see more details in Section~4) is a direct consequence of the theorem below.

\begin{theorem}
Let $A$ be a 3-free Artin group. Let $w$ be a geodesic word representing an element of $A$ and let $t$ be a letter. Then, there is a algorithm that takes $w$ and $t$ and returns a geodesic word equivalent to $wt$. {The complexity of this algorithm is linear on the length of $w$.}
\end{theorem}

This result will be formulated with more details in the Main Theorem of Section~4. The steps of our algorithm for solving the word problem are composed of free reductions and a generalization of the homogeneous relations (called $\tau$-moves) first defined by Holt and Rees in \cite{Holt2012}. As a corollary of this result, we obtain the following stronger result (which is stated with more detail in Corollary~\ref{Cor:LengthPreserved}):

\begin{corollary}
Let $g_1$ and~$g_2$ two geodesic words representing the same element of a 3-free Artin group. Then using only a finite set of homogeneous relations (called $\tau$-moves), we can transform~$g_1$ into~$g_2$. That is, we can transform~$g_1$ into~$g_2$ without increasing the length at any step.
\end{corollary}

Our work does not show that 3-free Artin groups are automatic. The strategies we employ to show that there is a solution to the word problem are similar to the ones that Holt and Rees use in \cite{Holt2012} and \cite{Holt2013} to show that large type and sufficiently large type Artin groups are automatic. Their work implies that geodesic words in large and sufficiently large Artin groups are automatic by proving a result similar to our \Cref{Cor:LengthPreserved}, and then showing that two geodesic words that are related by a single sequence of $\tau$-moves satisfy a fellow traveler property. Our methods generalize the concept of a $\tau$-move to include an arbitrary number of commutations in a single $\tau$-move, and so while we are still able to use these generalized $\tau$-moves to show that the word problem is solvable, we do not get the fellow traveler properties required for automaticity. 

 Another consequence of our work is that 3-free Artin groups satisfy Dehornoy's Property~H, a weak hyperbolicity property. Dehornoy originally defined groups satisfying Property H in 2011 \cite{Dehornoy2011} as groups whose word problem can be solved via an algorithm that is similar to Dehn's algorithm. In this sense, we can think of Property H as a weak hyperbolicity property \cite{DehornoyGodelle2013}. It is known that sufficiently large type Artin groups and FC-type Artin groups satisfy Property~H \cite{DehornoyGodelle2013, GodelleRees2016, DehornoyHoltRees2018} and conjectured that all Artin groups satisfy Property~H \cite{DehornoyGodelle2013}. Our work makes significant progress towards this conjecture.  (See Section~\ref{Sec:PropH} for a more precise definition of Property~H.) 

\begin{corollary} [Property~H]
Let $A_S$ be a 3-free Artin group. Then this group satisfies Property~H.
\end{corollary}

We are also aware of work in progress of Arya Juhasz that investigates the solvability of word problem for 3-free Artin groups using small cancellation theory and isoperimetric functions. This work differs from our approach in its methods and in that it does not include the construction of an explicit algorithm for determining if a word is geodesic as seen in our work. As of this writing, this work has not yet appeared on the ArXiv.

\medskip

This paper is subdivided in the following sections:  In Section 2, we give background and notation, and introduce the concept of a critical word. In Section 3, we define the concept of a Rightward Reducing Sequence, which will be the main component of the word problem solution algorithm. In Section 4, we prove that our Main Theorem depends on three main propositions. In Section 5, we prove these three propositions. Proofs of auxiliary lemmas are left to Section 6.

%% file: Background.tex
In this section, we recall some definitions and results from Holt and Rees in \cite{Holt2012}, \cite{Holt2013}. These articles study the word problem and automaticity in large and sufficiently large type Artin groups. The main premise of these papers is to partition a given word into subwords, each of which is contained in a dihedral Artin group. By applying facts about dihedral Artin groups to these subwords in a sequential fashion, Holt and Rees find the shortlex representative for an element of the Artin group. This section will cover the necessary background on dihedral Artin groups and define the notation needed going forward. Many of the results of Holt and Rees listed here have stronger versions given in their papers, but here we will only restate the results required for our work. 

\subsection{Definitions and notation for words in Artin groups}

\medskip 

We define an \textit{alphabet} to be a finite set $\mathcal{A}$. An element $a\in \mathcal{A}$ is called a \textit{letter}. A \textit{word} over $\mathcal{A}$ is a finite sequence of letters.
Formally, a word can be defined as a map $w : \{1,...,n\} \rightarrow \mathcal{A}$ where~$w (i)$ is the i-th letter of the word. The \textit{length} of a word~$w$ is the integer~$n$ and it is denoted by $|w|$. When $n=0$, we say that~$w$ is the \textit{empty} or \textit{trivial word} over~$\mathcal{A}$, and it is denoted by $\epsilon$. We denote by~$\mathcal{A}^{*}$ the set of all words over the alphabet~$\mathcal{A}$. 

\medskip

Given a word $w=w_1w_2w_3$, with possibly empty words $w_1, w_2, w_3\in \mathcal{A}^{*}$, the word~$w_1$ is said to be a \textit{prefix} of~$w$, $w_3$ is a \textit{suffix} of~$w$ and~$w_2$ is a \textit{subword} of~$w$. Given a word~$w$ we denote by~$f[w]$ and~$l[w]$ the first and last letters of $w$ respectively. 
We also denote $\mathrm{suf}[w]$ the word $w$ without its first letter and $\mathrm{pref}[w]$ the word~$w$ without its last letter.

From now on we fix $\mathcal{A}=S\cup S^{-1}$, where~$S$ is the standard generating set of the Artin group~$G$. A letter $a\in \mathcal{A}$ is \textit{positive} if $a\in S$ and is \textit{negative} otherwise. The \textit{name} of a letter is its positive form. We say that two letters are \textit{essentially different} if they have different names. 
We say that a word $w\in \mathcal{A}^{*}$ is \textit{positive} if all its letters are positive, \textit{negative} if all its letters are negative and \textit{unsigned} otherwise.

\medskip

A word $w\in \mathcal{A}^{*}$ is \textit{freely reduced} if it does not admit any subword of the form $aa^{-1}$ or $a^{-1}a$ for any letter~$a$. We say that a word that is not freely reduced word \textit{admits a free reduction}.
	
\medskip	
	
Two words $v$ and $w$ represent the same element of the Artin group if we can change $v$ into $w$ via a sequence of insertions of words of the form $aa^{-1}$ or $a^{-1} a$, free reductions, or applications of the Artin relations. If $v$ and $w$ represent the same group element, we denote this by $w =_{G} v$ and say that~$w$ and~$v$ are \emph{equivalent in~$G$}. A word $w\in \mathcal{A}^{*}$ is \textit{geodesic} if for any other word~$v$ such that $w =_{G} v$, we have that $|w|\leq |v|$.

\medskip

We say that a letter~$a$ \textit{appears in} the word~$w$, if either~$a$ or~$a^{-1}$ is a letter in~$w$. If there are letters $a,b,c$ in a word~$w$ such that~$a$ appears before~$b$,~$b$ appears before~$c$ and~$b$ does not commute with either ~$a$ or ~$c$ we say that~$b$ is \emph{trapped} in~$w$, {because we cannot move~$b$ out from between the letters~$a$ and~$c$ using only commutations between individual letters}. 

\medskip

\subsection{Geodesics in dihedral Artin groups}

Following the notation of Holt and Rees \cite{Holt2012, Holt2013}, we define the notation ${}_{m}(a,b)$ as the alternating product of the letters~$a$ and~$b$ of length~$m$, starting with the letter~$a$. Similarly, we define $(a,b)_{m}$ as the alternating product of~$a$ and~$b$ of length~$m$ ending with $b$. Notice that if~$m$ is even, we have $_{m}(a,b)=(a,b)_{m}$ and we can use the two expressions interchangeably.

\begin{definition}
	
	The \textit{dihedral Artin group}  (also called the \textit{two-generator Artin group}) $A_{2}(m)$, $m\in\{ 2,3,\dots +\infty\}$ is the Artin group having two generators with a relation of length $m$ between them if $m\in\mathbb{Z}$ or the free group with two generators otherwise. If $m<\infty$ this is the group with presentation $\langle a,b \mid {}_{m}(a,b)={}_{m}(b,a)\rangle$.
	
\end{definition}

We say that the words that contains only instances of two essentially different letters are \emph{2-generated}. Words that are 2-generated represent elements of a dihedral Artin group. 

\medskip

In \cite{Mairesse2006}, Mairesse and Matheus give a method for identifying geodesic words in dihedral Artin groups. Let $w$ be a freely reduced word in $A_{2}(m)$ over the alphabet $\mathcal{A}=\{a, a^{-1}, b, b^{-1}\}$. Consider the integers:
$$	
\begin{array}{rcl}
r_{1} &=&\max\{r\mid {}_{r}(a,b) \text{ or } _{r}(b,a) \text{ is a subword of } w \} ;     \\
r_{2}&=&\max\{r\mid {}_{r}(a^{-1},b ^{-1}) \text{ or } _{r}(b^{-1},a^{-1}) \text{ is a subword of } w \}; \\
p(w)&=&\min\{r_{1},m\}; \\ n(w)&=&\min\{r_{2},m\}.
\end{array}
$$
The following result proves that geodesic words~$w$ in~$A_{2}(m)$ are characterized by the values~$p(w)$ and~$n(w)$:

\begin{proposition}\label{prop:pn}\cite{Mairesse2006} Let $g\in A_{2}(m)$ and let $w\in \mathcal{A}^{*}$ be a freely reduced word representing $g$.
	
	\begin{enumerate}
		\item If $p(w)+n(w)<m$, then $w$ is the unique geodesic representative for $g$.
		\item If $p(w)+n(w)=m$, then $w$ is one of the geodesic representatives for $g$.
		\item If $p(w)+n(w)>m$, then $w$ is not geodesic. Furthermore, $w$ has a subword $w'$ such that $p(w')+n(w')=m$.
		
	\end{enumerate}
\end{proposition}

We are particularly interested in words which have more than one geodesic representative. Following \cite{Holt2013}, we define a special subset of these words, called critical words. 

\begin{definition}\label{def: critical words}\cite{Holt2012}
	Let $w$ be a freely reduced word in $A_{2}(m)$. Let $\{x,y\}=\{z,t\}$ be sets of generators of $A_2(m)$ and write $p=p(w)$ and $n=n(w)$. The word $w$ is called a \textit{critical word} if $p+n=m$ and it has the form
	\[w= {}_{p}(x,y)\eta(z^{-1},t^{-1})_{n} \text{ or } w= {}_{n}(x^{-1},y^{-1})\eta(z,t)_{p}\]	
	
	Here $\eta$  represents any word in~$\mathcal{A}^{*}$. Notice that~$p$ or~$n$ could be~$0$. In this case we have that $(a,b)_0=\,_0(a,b)=\epsilon$. Moreover, if~$w$ is negative or positive, we add the restriction that~$w$ can only contain one subword of the form $(x,y)_m$ or $(x^{-1},y^{-1})_m$.

\end{definition}

Note that while the definition does not explicitly give restrictions on~$\eta$, the requirement that $p+n=m$ implicitly imposes restrictions on~$\eta$. We will use~$\eta^+$ (resp. $\eta^-$) to express that~$\eta^+$ is positive (resp. negative).

\begin{lemma}\label{lemma: critical_subwords}
If $w$ is a 2-generated critical word that contains a critical subword $\zeta$, then the smallest suffix of $w$ containing $\zeta$ has a critical suffix, $\xi$, with $p(w)=p(\xi)$ and $n(w)=n(\xi)$. 
\end{lemma}

\begin{proof}
Assume that $w$ is positive i.e. $p(w)=m$ (the negative case is analogous). The subword $\zeta$ must be positive and must contain the only instance of the alternating positive word of length $m$. So $\zeta$ is a prefix of $w=(a,b)_m \eta^+$ or a suffix of $w=\eta^+ (a,b)_m$ and the result follows.

Assume that $w$ is neither positive nor negative. We know that $n(\zeta)=n(w)$ and $p(\zeta)=p(w)$. Suppose that $w={}_{p}(x,y)\eta(z^{-1},t^{-1})_{n}$ (the other case is symmetrical).  In this case the desired critical suffix of $w$ starts at the first letter of a positive alternating word of length $p(w)$ inside $\zeta$. 
\end{proof}

We now define an involution $\tau$ on the set of critical words, again following the work of Holt and Rees \cite{Holt2012} and Brien \cite{Brien2012}. For a dihedral Artin group generated by $a$ and $b$, consider the homomorphism $\delta$ defined on the set of words as follows: if $m$ is even then $\delta(a)=a, \delta(b)=b$, if $m$ is odd then $\delta(a)=b, \delta(b)=a$.

Then we define, for $n,p>0$:

$$\tau\left({}_{p}(x,y)\eta(z^{-1},t^{-1})_{n}\right)= \,_n(y^{-1},x^{-1})\delta(\eta)(t,z)_p,$$ $$
\tau\left({}_{p}(x^{-1},y^{-1})\eta(z,t)_{n}\right)= \,_n(y,x)\delta(\eta)(t^{-1},z^{-1})_p.$$

For the cases $p=0$ and $n=0$, the critical words are negative or positive and we define:

$$\tau\left(\eta^-(z^{-1},t^{-1})_{m}\right)= \,_m(y^{-1},x^{-1})\delta(\eta^{-}), \quad y^{-1}\neq f[\eta^-],$$
$$\tau\left({}_{m}(x^{-1},y^{-1})\eta^{-}\right)= \delta(\eta^-)(t^{-1},z^{-1})_m, \quad z^{-1} \neq l[\eta^-],$$

\smallskip
$$\tau\left(\eta^+(z,t)_{m}\right)= \,_m(y,x)\delta(\eta^{+}), \quad y\neq f[\eta^+],$$
$$\tau\left({}_{m}(x,y)\eta^{+}\right)= \delta(\eta^+)(t,z)_m, \quad z \neq l[\eta^+].$$

{For example, if $m=5$ we have }

{$$\underbrace{ababa}_{{}_m(a,b)}\underbrace{a^3ba}_{\eta^+}\xmapsto{\tau} \underbrace{b^3ab}_{\delta(\eta^+)}\underbrace{babab}_{(a,b)_m}$$
\smallskip
$$\underbrace{aba}_{{}_p(a,b)}\underbrace{a^3b^{-2}a}_{\eta}\underbrace{b^{-1}a^{-1}}_{(b^{-1},a^{-1})_n}\xmapsto{\tau} \underbrace{b^{-1}a^{-1}}_{{}_n(b^{-1},a^{-1})} \underbrace{b^3a^{-2}b}_{\delta(\eta)}\underbrace{bab}_{(a,b)_p}$$}

The technical definition of $\tau$ will be used just in a few proofs. For the rest of them, the reader should focus on the important properties of this map given below. 

\begin{proposition}\label{prop:tau properties} \cite{Holt2012,Brien2012} $\tau$ satisfies the following properties for any critical word $w$:
	
	\begin{enumerate}
		\item $\tau(w)$ and $w$ represent the same group element, $\tau(w)$ is a critical word and $\tau(\tau(w))=w$.
		
		\item $p(\tau(w))=p(w)$ and $n(\tau(w))=n(w)$.
		
		\item $f[w]$ and $f[\tau(w)]$ have different names, the same is true for $l[w]$ and $l[\tau(w)]$.
		
		\item $f[w]$ and $f[\tau(w)]$ have the same sign if $w$ is positive or negative, but different signs if $w$ is unsigned; the same is true for $l[w]$ and $l[\tau(w)]$.
	\end{enumerate}
	
\end{proposition}

Given a word $w$ representing an element in a dihedral Artin group with $u$ a critical subword, we can replace $u$ with $\tau(u)$ without changing the length of the word or the represented group element. We call such a replacement a \emph{$\tau$-move}. Holt and Rees show that any non-geodesic word in a dihedral Artin group can be reduced to a geodesic representative by a sequence of free reductions and $\tau$-moves. They then generalize this method to large and sufficiently large Artin groups. We plan to use this method to study 3-free Artin groups. 

\medskip

We finish this section by proving the following lemma, which will be an important step in our arguments going forward. 

\begin{lemma}\label{lemma: critical_subwords2}
Let~$w$ be a 2-generated critical word that contains a critical subword~$\zeta$. Let~$w'$ be the word obtained from~$w$ by replacing~$\zeta$ with $\tau(\zeta)$. Then the smallest suffix of $w'$ containing~$\tau(\zeta)$ contains a critical suffix~$\xi$ with $p(w)=p(\xi)$ and $n(w)=n(\xi)$.
\end{lemma}

\begin{proof}
Assume that $w$ is positive (the negative case is analogous). The subword $\zeta$ must be a positive critical subword and $p(\zeta)=p(\tau(\zeta))=m$. Thus there is at least one positive alternating sequence of length~$m$ in~$w'$. Let~$\xi$ be the suffix of~$w'$ starting at the rightmost such sequence. The suffix~$\xi$ will be a positive word generated by~$a,b$ with exactly one positive alternating sequence of length~$m$ an hence a critical word.

Assume that~$w$ is unsigned. 
We know that $p(\zeta)=p(w)$ and $n(\zeta)=n(w)$ so~$\zeta$ is also an unsigned critical word and~$\zeta$ must begin in an alternation. Any letter in~$w$ before this alternation should not increase the length of the alternation. Thus, for example if $f[\zeta]=a$, then either $\zeta$ is a prefix of $w$ or the letter preceding $a$ in~$w$ must be in $\{a,b^{-1}\}$. Applying~$\tau$ to~$\zeta$ we get $f[\tau(\zeta)]=b^{-1}$. Thus the alternation at the beginning of $\tau(\zeta)$ does not introduce any free cancellations to $w'$ nor does it increase the value of $p$ or $n$.  We also know that $\zeta$ ends in an alternation and a similar argument can be applied to the alternation at the end of $\tau(\zeta)$. Thus~$w'$ is also a freely reduced word with the same value of~$p$ and~$n$.

Suppose that $w={}_{p}(x,y)\eta(z^{-1},t^{-1})_{n}$ (the other case is symmetrical). Because there is a positive alternation of length $p$ in $w$, there must also be a positive  alternation of length $p$ in $w'$. Let $\xi$ be a suffix starting at such an alternation. Then $\xi$ will be a freely reduced word with $p(\xi)=p(w), n(\xi)=n(w)$ and such that $\xi$ begins with a positive alternation of length $p$ and ends with a negative alternation of length $n$. Thus $\xi$ is the desired critical suffix. 
\end{proof}

%% file: RRS.tex
The main aim of this section is to generalize our discussion to 3-free Artin groups, via sequences of transformations of subwords called \emph{rightward reducing sequences}. For an Artin group $A_S$ we denote $\mathcal{A}^*$ the set of words in $\mathcal{A}=S\cup S^{-1}$.

\medskip

We wish to use our knowledge of dihedral Artin groups, so we will now define a way of using commutations to isolate a 2-generated subword inside a larger word. 
Informally, we will say that a word $w\in \mathcal{A}^*$ is \emph{pseudo 2-generated} (P2G) in the pseudo-generators $a$ and $b$ if, for every instance of a letter in $w$ not in $\{a,b,a^{-1},b^{-1}\}$, we can use commuting relations to push this letter either to right or to the left of all the letters in $\{a,b,a^{-1},b^{-1}\}$. This effectively isolates the pseudo-generators in the center of the word and allows us to treat this central subword as a word in a dihedral Artin group. 

\medskip

Recall that, given a letter, the name of a letter is its positive form and we say that two letters are essentially different if they have different names.  

\medskip

The formal definition of a P2G word is as follows.

\begin{definition} 
	Let $a,b\in S$ with $2<m_{a,b}<\infty$ and denote by $P$ the set $\{a,b,a^{-1},b^{-1}\}$. Let $w\in \mathcal{A}^*$ be a word such that $f[w]$, $l[w]\in P$. 
	Let $w_p$ be the prefix of $w$ up to but not including the first instance of a letter in $P$ that is essentially different from $f[w]$. Similarly, let $w_s$ be the largest suffix of $w$ which does not contain a letter in $P$ that is essentially different from $l[w]$.  Factor~$w$ as ${w}=w_pw_qw_s$.

	We say that $w$ is a \textit{pseudo $2$-generator word} (P2G word in the following) in \textit{pseudo-generators} $\{a,b\}$ if all the letters in~$w_p$ commute with~$f[w]$, all the letters in~$w_q$ not in~$P$ commute with both~$a$ and~$b$, and all the letters in~$w_s$ commute with~$l[w]$. The letters of $w$ not in $\{a,b,a^{-1},b^{-1}\}$ are called \textit{internal letters}.
	
	\end{definition}	
	
\begin{example} Suppose that $a,b\in S$ such that $m_{a,b}=4$. Further suppose that $x,y,z\in S$ such that $x$ commutes with $a$, $y$ commutes with $b$, and $z$ commutes with both $a$ and $b$. Then the word $(axzax)(bza)(byb^{-1})$ is P2G  with $w_p=axzax$, $w_q=bza$ and $w_s=byb^{-1}$. Notice that there can also be P2G words that begin and end in the same generator, for example $(axa)( bzab )(axa)$ is also P2G with~$w_p$, $w_q$, and $w_s$ indicated by parentheses. 
\end{example}
	
\medskip	

 
 We now define some further notation associated to a given P2G word $w$. If $w$ is a P2G word, as above, then it is equivalent to a word $\alpha\rho\hat{w}\beta$ where: 
 
 \begin{itemize}
  \item $\hat{w}$ is the word obtained from $w$ after deleting all the internal letters.
	\item $\alpha$ is the word obtained from $w_p$ after deleting all instances of pseudo-generators. By definition of a P2G word, all the letters in $\alpha$ must commute with $f[w]$ and may or may not commute with the other pseudo-generator. 
	\item $\beta$ is the word obtained from $w_s$ by deleting all the pseudo-generators and all the letters that can be pushed to the left of $\hat{w}$ in $w$ via commutations. All the letters in $\beta$ must be able to be pushed to the right via commutations and so must commute with $l[w]$.
	\item $\rho$ is the word obtained from $w$ by deleting all the pseudo-generators and all the letters in $\alpha$ and $\beta$. By the definition of $\alpha$, $\beta$ and P2G words, all the letters in $\rho$ must commute with both pseudo-generators.  
		
	\end{itemize}
	
	The word $\alpha\rho\hat{w}\beta$ can be obtained from $w$ by pushing as many internal letters as possible to the right via commuting relations, and pushing the remaining internal letters to the left via commuting relations. 

%
	
		

When dealing with several P2G words, we will use subindices.

\begin{example}
 Suppose that $a,b\in S$ such that $m_{a,b}=4$. Further suppose that $x,y,z\in S$ such that $x$ commutes with $a$, $y$ commutes with $b$, and $z$ commutes with both $a$ and $b$ and with $y$.	The word $w=(axzax)(bza)(byzb)$ is equivalent to~$\alpha_w\rho_w\hat{w}\beta_w$ where $\alpha_w=xzx,\, \rho_w= z^2,\, \hat{w}=a^2bab^2$ and $\beta_w=y$.
\end{example}


\medskip	
	
	We are especially interested in P2G words where the central subword $\hat{w}$ is a critical word in the sense of Holt-Rees. This allows us to generalize the notion of a critical word beyond dihedral Artin groups. 
	
	\begin{definition}
	Let $w$ be a P2G word in $\{a,b\}$. Let $\hat{w}$ be the word obtained from $w$ by erasing the internal letters. We say that $w$ is \textit{P2G critical} if $\hat{w}$ is a 2-generated critical word in the sense of Holt and Rees. See \Cref{def: critical words} above. 
\end{definition}

\begin{remark}
Notice that a P2G critical word, $w$, entirely determines its  pseudo-generators: the pseudo-generators will always be the first letter of $w$, and the first letter after $f[w]$ that does not commute with $f[w]$.
\end{remark}

Now we extend the involution $\tau$ on 2-generated critical words defined in Section~\ref{Section: Background} so that it will apply to any P2G critical word. 

\begin{definition}
	Given a P2G critical word $w$ as above we define the $\tau$-move:
	$$\tau(w)=\alpha \rho \tau(\hat{w}) \beta,$$
	where $\alpha, \rho, \hat{w}, \beta$ are defined as above and $\tau(\hat{w})$ is the $\tau$-move defined for 2-generated critical words by Holt-Rees in \cite{Holt2012} explained in Section \ref{Section: Background}. We say that this $\tau$-move \emph{produces} the letter $l[\tau(\hat{w})]$. 
	
	
\end{definition}

Notice that given a P2G critical word $w$, the properties of $\tau(w)$ listed in \Cref{prop:tau properties} may not apply. 
However, the following corollary of \Cref{prop:tau properties} (proven in \cite{Holt2012} for 2-generated critical words) extends to the P2G critical case. 

\begin{corollary}\label{Cor:sufixes_same_letters}
If $w$ is a P2G critical word and $w'$ is a prefix of $w$ which is also P2G critical in the same pseudo-generators, then $f[\tau(w)]=f[\tau(w')]$. Similarly, if $w$ is a P2G word and $w'$ is a P2G critical suffix in the same generators, then $\beta_w=\beta_w'$ and $l[\tau(\hat{w})]=l[\tau(\hat{w}')]$.
\end{corollary}

\begin{proof}
This follows directly from  \Cref{prop:tau properties},  \cite[Corollary 2.2]{Holt2012} and the definition of a P2G word.
\end{proof}

One particularly nice feature of $\tau$-moves, in both the dihedral case and when applied more generally to P2G critical words, is that $w$ and $\tau(w)$ are words of the same length. We would like to find an algorithm which begins with a word $w$ and reduces this word to a geodesic by only applying free reductions and $\tau$-moves to subwords of $w$. This will give us a way of obtaining a geodesic representative of $w$ without ever increasing the length at any intermediate step. 
 To that end, we now define a way of applying $\tau$-moves to a word in sequence. This sequence of $\tau$-moves will take a non-geodesic but freely reduced word and apply successive $\tau$-moves until a free reduction becomes possible. 
  
\begin{figure}[ht]
\centering
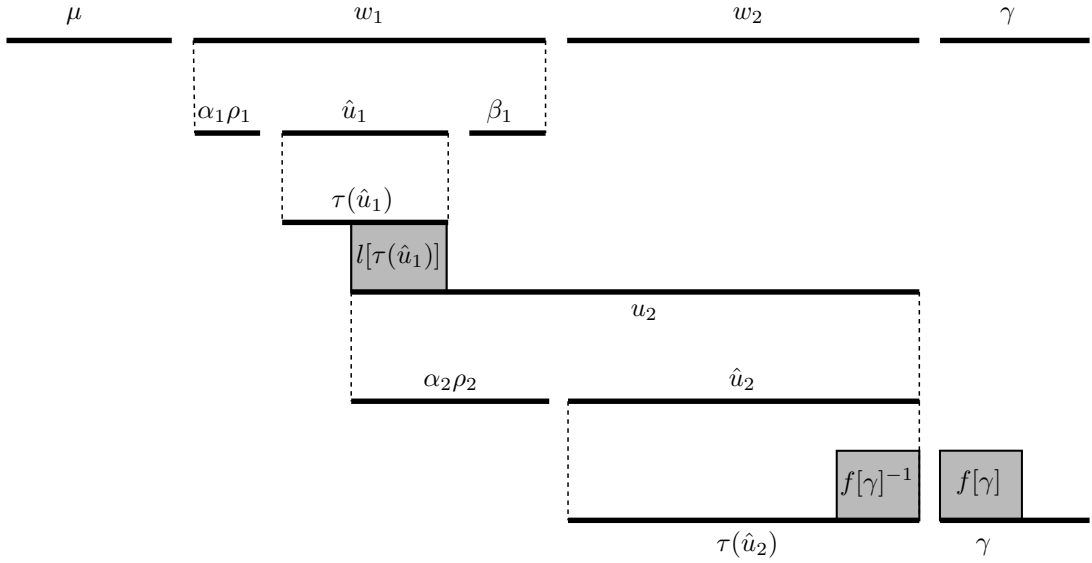
\caption{Diagram for an RRS $U=u_1,u_2,u_3$ with $w_3=\epsilon$. (Notice by Lemma~\ref{lemma:beta} that $\beta_2=\epsilon$).}\label{RRS1}
\end{figure}

\begin{figure}[ht]
\centering
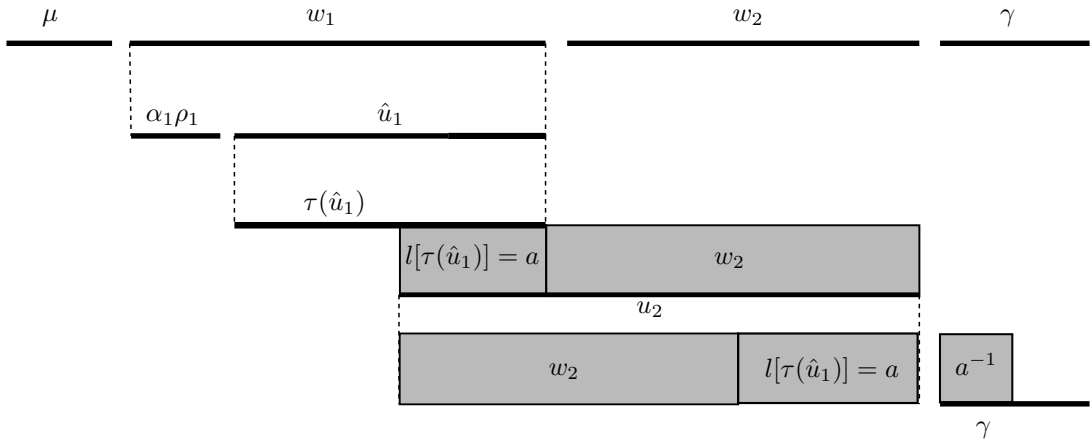
\caption{Diagram for an RRS $U=u_1,u_2$ with $u_2=aw_2$. (Notice by Lemma~\ref{lemma:beta} that $\beta_1=\epsilon$).}\label{RRS2}
\end{figure}

\begin{definition}[Rightward Reducing Sequence] 
The reader can find pictures for this definition in Figures~\ref{RRS1} and \ref{RRS2}. Let $w$ be a freely reduced word and $U=u_1,\dots,u_k,u_{k+1}$ a sequence of words such that, for $i\leq k$, $u_i$ is a P2G critical word, and $u_{k+1}=av$, where $a$ is a letter commuting with every letter in the word $v$. 

We note $\tau(u_i)= \alpha_i\rho_i\tau(\hat{u}_i)\beta_i$ for all $i\leq k$ and in particular $l[\tau(\hat{u}_i)]\beta_i$ is a suffix of $\tau(u_i)$. We say that $w$ admits $U$ as a \textit{rightward reducing sequence} (RRS) of length $k$ if $w$ can be written as a product of words $w=\mu w_1 \dots w_k w_{k+1}\gamma$ such that $w_i\neq \epsilon$ for $i\leq k$, $w_1=u_1$, $u_i=l[\tau(\hat{u}_{i-1})]\beta_{i-1}w_i$ for $1<i\leq k+1$, and $f[u_{k+1}]=f[\gamma]^{-1}$. Notice that $w_{k+1}$ can be trivial and that if $k=0$ then $w_1=f[\gamma]^{-1} v$.

We respectively call $\mu$ and $\gamma$  the \emph{head} and the \emph{tail} of the RRS.

\end{definition}

\medskip

In this definition, we begin with a non-geodesic word $w=\mu w_1 \dots w_kw_{k+1}\gamma$. We assume that $w_1$ is P2G critical and apply a $\tau$-move to $w_1$. This produces a new word $\mu \tau(w_1)w_2\dots w_kw_{k+1}\gamma$. We now have a P2G critical subword $u_2=l[\tau(\hat{u}_1)]\beta_1w_2$ and we apply $\tau$ to this subword. We repeat this step, applying successive $\tau$-moves and maintaining the length of the word at each step. After the final $\tau$-move we can transform $u_{k+1}=aw_{k+1}$ in $w_{k+1}a$ and obtain a word which can be freely reduced. 

Intuitively, we can think of an RRS as a way of applying successive $\tau$-moves and commutations, moving left to right, until a free reduction becomes possible. Notice that any word which is not freely reduced admits an RRS with $k=0$ and $w_1$ a cancelling letter. In this case, the RRS makes no changes to the word. 
 
 \begin{example}
 
 Suppose all the $u_i$ have the property that~$f[u_i]$ and~$l[u_i]$ are essentially different. Then we can suppose that $a_i$ and $b_i$ are the pseudo-generators of~$u_i$ with~$a_i$ the name of~$f[u_i]$ and~$b_i$ the name of~$l[u_i]$. We know that $\tau$-moves on two-generator words always change the name of the last letter of the critical word (Proposition~\ref{prop:tau properties}), so the name of $l[\tau(\hat{u}_i)]=f[u_{i+1}]$ is $a_i$. Thus $a_i=a_{i+1}$. In this special case, the pseudo-generator $a_1$ is the name of the first letter in the first critical word in the RRS. Each successive $\tau$-move pushes this letter further to the right in the word~$w$, possibly changing its sign in the process. Each $\tau$-move  also makes other changes to the word $w$, however these changes do not affect the length of~$w$.  After the final $\tau$-move and (possibly) some commutations, the new word will contain a free reduction in~$a_1$. Thus, in this special case, the RRS is effectively moving a single letter to right, (while making other changes at the same time) until a cancellation is possible. This means that some RRS effectively use $\tau$-moves to move a single letter further to the right in the word.

For another more concrete example, consider an Artin group~$A$ with four generators $a,b,c,d$ where~$a$ and~$b$ have a relation of length~4, $b$ and~$c$ a relation of length~5, $d$ commutes with~$a$, $b$ and~$c$, and~$c$ commutes with~$a$. The word $w=(acbab^2cda)(cb^{-1}c^{-1}b^{-1})d^5 c^{-1}$ admits an RRS $U=u_1,u_2,u_3$ with factorization $w=\mu w_1w_2w_3 \gamma$ where $\mu=\epsilon$, $w_1=acbab^2cda$, $w_2=cb^{-1}c^{-1}b^{-1}$, $w_3=d^5$ and $\gamma=c^{-1}$. By definition, $w_1=u_1$ so $u_1=_A \alpha_1\rho_1 \hat{u}_1 \beta_1 $ with $\alpha_1=c$, $\rho_1=d$, $\hat{u}_1=abab^2a$ and $\beta_1=c$. When we apply~$\tau$ we get $\tau(\hat{u}_1)=ba^2bab$, hence $u_2=bc^2b^{-1}c^{-1}b^{-1}= \hat{u}_2$. We apply~$\tau$ again and have  $\tau(bc^2b^{-1}c^{-1}b^{-1})=c^{-1}b^{-1}c^{-1}b^2c$, so $u_3=cd^5$. Finally, $u_3$ can be changed via commutations to $d^5c$. After the RRS we obtain a new word $cdba^2bac^{-1}b^{-1}c^{-1}b^2d^5cc^{-1}$ and the last two letters give a free cancellation. 

 \end{example}


\smallskip

Recall that, if $w$ is  P2G word in pseudo-generators $a$ and $b$, then $\hat{w}$ is the word obtained from $w$ by erasing internal letters in $w$. The next lemma shows that every P2G word $w$ with non geodesic $\hat{w}$ admits an RRS:

\begin{lemma}\label{Lemma: non-critical P2G has RRS}
Suppose that $w$ is a freely reduced P2G word and $\hat{w}$ is not geodesic in the dihedral Artin group~$A_{\{a,b\}}$.  Let $\hat{v}$ be the largest prefix of $\hat{w}$ such that $\hat{v}$ is geodesic, and consider the nongeodesic prefix of $\hat{w}$ given by $\hat{v}x$ where $x\in \{a,a^{-1},b,b^{-1}\}$. 
Then $w$ admits an RRS of length 1, 
with  $l[\tau(\hat{w}_1)]=x^{-1}$.
\end{lemma}

\begin{proof}
Let $a$ and $b$ be the pseudo-generators of $w$. Notice that any subword of $w$ which begins and ends  in the letters in  $\{a,a^{-1},b,b^{-1}\}$ must be a P2G word in the same pseudo-generators. Similarly, given a subword $\hat{u}$ in $\hat{w}$, we can consider the corresponding subword $u$ in $w$ such that erasing all the internal letters from $u$ produces $\hat{u}$. 

By \cite[Lemma 2.3]{Holt2012}, there exists some suffix of $\hat{v}$, denoted by $\hat{w}_1$ such that $\hat{w}_1$ is a critical word and $l[\tau(\hat{w}_1)]=x^{-1} $. {The suffix $\hat{w}_1$ has a corresponding P2G critical subword of $w$ denoted by $w_1$. Let $w_2$ be the subword of $w$ beginning immediately after the last letter of $w_1$ and ending at $x$. Because of the way we have defined $w_1$ the word $w_2$ must consist entirely of internal letters in $w$. Because $w$ is P2G and $x$ or $x^{-1}$ appears in $w_1$, any letters in $w_2$ must commute with $x$. Otherwise they are internal letters trapped between pseudo-generators. This means that $w_1,w_2$ are subwords of $w$ that produce an RRS of length 1, which cancels the letter $x$.}
\end{proof}

The following lemma is a key characteristic of RRSs that is used repeatedly throughout the subsequent proofs and reveals one important way that the 3-free condition contributes to the main results of this paper. 

We denote by $t^\#$ some arbitrary nonzero power of the letter $t$.

\begin{lemma}\label{lemma:beta} Let $A_S$ be an Artin group, and Let $U = u_1,\dots, u_{k+1}$ be an RRS of a word $w$ with factorization $w=\mu w_1\cdots w_{k+1}\gamma$. Then:  
\begin{enumerate}
\item $\beta_k = \epsilon$.

\end{enumerate}

If we further suppose that $u_i$ and $u_{i+1}$, $1\leq i< k$, are two consecutive P2G critical words in U, and assume that $\beta_i\neq \epsilon$, then 
\begin{enumerate}
\setcounter{enumi}{1}
\item $\alpha_{i+1}=\epsilon$;
\item If $z,t$ are the pseudo-generators of $u_{i+1}$, and $m_{z,t}>3$, then $\beta_i=t^\#$ or $\beta_i=z^\#$.

\end{enumerate}

\end{lemma}

\begin{proof}

For the first item, notice {that by the definition of a P2G word $u$ and $\beta_u$, the first letter of $\beta_u$ must commute with $l[u]$, but not with the other pseudo-generator of $u$, which is $l[\tau(\hat{u})]$.} Thus if $\beta_k$ is not trivial, then {the first} letter in $\beta_k$ would lie in $u_{k+1}$ but would not commute with $l[\tau(\hat{u}_{k})]$, producing a contradiction.

\medskip

For the last two items, let $z_i$ and $t_i$ be the pseudo-generators of~$u_i$. We know that $l[\tau(\hat{u}_i)]=f[u_{i+1}]$ is a pseudo-generator of both~$u_i$ and~$u_{i+1}$, so we can assume that $l[\tau(\hat{u}_i)]=z_i^{\pm{1}}=z_{i+1}^{\pm{1}}$. We also know from the definition of P2G critical word that~$f[\beta_i]$ does not commute with~$z_{i+1}$. This means that when we split~$u_{i+1}$ into three blocks, as in the definition of P2G words, we see that the first block~$w_p$ consists only of the letter $z_{i+1}^{\pm 1}$ and so $\alpha_{i+1}=\epsilon$. 

\medskip

We can now write $u_{i+1}=z_{i+1}^{\pm 1} t_{i+1}^\#yw_{i+1}$, where $y$ is the suffix of the word $\beta_i$ after $t_{i+1}^\#$ such that $f[y]\neq t_{i+1}^{\pm 1} $. Note that $z_{i+1}$ cannot appear as a letter in $y$, as $z_{i+1}=z_{i}$ is the generator of $u_i$ which does not appear as a letter in $\beta_i$. However, since $\hat{u}_{i+1}$ must to be critical, $z_{i+1}$ must appear at least once in $yw_{i+1}$, so $z_{i+1}$ is a letter in $w_{i+1}$. The fact that the Artin relation between $t_{i+1}$ and $z_{i+1}$ is of length greater than 3 implies that there must also be at least one instance of $t_{i+1}$ in $w_{i+1}$.

In order for $u_{i+1}$ to be a P2G word with pseudo-generators $z_{i+1}$ and $t_{i+1}$, we must be able to push any letters in $y$ not equal to $t_{i+1}^{\pm 1}$ either to the left or the right using commutations. However, by definition of $\beta_i$, $f[y]$ cannot be pushed to the left of $\hat{u}_i$, so either $f[y]$ does not commute with $t_{i+1}$ or $f[y]$ does not commute with $z_i=z_{i+1}$. This means that~$f[y]$ is trapped in $u_{i+1}$ and in order for $u_{i+1}$ to be a P2G critical word, $y=\epsilon$ as desired. 
\end{proof} 
 
 \medskip

\subsection{Optimal Rightward Reducing Sequence}
 
As one may apply different RRSs to the same word, we will determine an order of application of RRSs by using the following definition:

\begin{definition}
	Let $U=u_1,\dots,u_k,u_{k+1}$ be an RRS of a word $w$ with associated decomposition $w=\mu w_1\dots w_k w_{k+1}\gamma$. Let $z_l$ and $t_l$ be the two pseudo-generators of~${u}_l$ for $1\leq l \leq k$. We call an RRS \textbf{optimal} if the following conditions are satisfied:
	
	\begin{enumerate}
		\item  		
	{The first letter of $w_1$ is at least as far to the right as in any other RSS of~$w$}. In other words, $\mu$ is as long as possible.
	
	\item  	
	{ $l[\tau(u_i)]\neq f[w_{i+1}]^{-1}$ for $i\leq k$ and $f[\gamma]$ does not appear in $w_{k+1}$}. In other words, the RRS does not introduce any free reduction, except at $l[\tau(u_k)]=f[\gamma]^{-1}$. 

		\item If 
		{every letter in $\alpha_l$ commutes with both $z_l$ and $t_l$} then $\mid\{{z_{{l-1}},t_{{l-1}}\}\cap \{z_{{l}},t_{{l}}\}\mid=1}$.
	\end{enumerate} 
\end{definition}
 
 \noindent
 Intuitively, condition 1) guarantees that the head $\mu$ is as long as possible, condition 2) guarantees that the tail $\gamma$ is as long as possible, and condition 3) guarantees that the length $k$ is minimized.
 
\begin{remark}\label{remark:optimality}
If $\beta_i=\epsilon$ and $l[\tau(\hat{u}_i)]$ has the same name as the first pseudo-generator in~$w_{i+1}$ then the RRS is not optimal. This follows from the fact that optimality implies that the critical word~$\hat{u}_{i+1}$ begins in an alternation that contributes to~$p(\hat{u}_{i+1})$ or~$n(\hat{u}_{i+1})$.
\end{remark}

 \begin{lemma}\label{ExistenceOfOptimalRRS}
For a 3-free Artin group, if $w$ admits an RRS, then it admits a unique optimal RRS.
 \end{lemma}

\begin{proof}
We first prove existence.
Suppose that the word $w$ admits an RRS $U=u_1,\dots,u_{k+1}$ with a decomposition $w=\mu w_1\cdots w_k w_{k+1} \gamma$. We can suppose that this RRS is the one with the head $\mu$ as long as possible. If there is some $i$ such that $l[\tau(u_i)]=f[w_{i+1}]^{-1}$, redefine  $u_{i+1}=l[\tau(u_i)]$, and obtain a new RRS $U=u_1,\dots u_i, u_{i+1}$ with the decomposition $w=\mu w_1\cdots w_i\gamma'$. Here $\gamma'$ is the suffix of $w$ starting at the first letter of $w_{i+1}$ from the original RRS and the new RRS has $w_{i+1}=\epsilon$.  Finally, if $f[\gamma]$ appears in $w_{k+1}$, redefine $w_{k+1}$ to end at just before this instance $f[\gamma]$. We now have an RRS that satisfies conditions 1 and 2 of optimality. 

Suppose that the third optimality condition is not satisfied. Consider an integer $l$ such that $z,t$ are the pseudo-generators for both~$u_{l-1}$ and~$u_l$, and suppose that every letter in~$\alpha_l$ commutes with both~$z$ and~$t$. 
Consider the word $v=u_{l-1}w_{l}$. This must also be a P2G word with pseudo-generators $z,t$: the assumptions about~$l$ imply that that~$\beta_{l-1}$ is trivial (Lemma~\ref{lemma:beta}) and all letters not equal to~$z$ and $t$ can escape via commutations. We divide into cases according to whether or not $\hat{v}$ is geodesic. 

If $\hat{v}$ is not geodesic, then \Cref{Lemma: non-critical P2G has RRS} implies that there is a proper prefix~$v'$ of~$\hat{v}$ that is critical and produces {an RRS of length 1 that cancels the first instance of $t$ or $z$} immediately to the right of~$v'$. Redefine $\hat{u}_{l-1}=v'$ {and redefine $u_l$ to consist of the letters in $v$ between $v'$ and the cancelled letter}. Finally, redefine the tail~$\gamma$ to start at the previously mentioned cancelled letter.

If $\hat{v}$ is geodesic, then we will show that the word $\hat{v}$ is also a critical word, that $\beta_v=\beta_{u_l}$, and $l[\tau(\hat{v})]=l[\tau(\hat{u}_l)]$. This means we can replace the RRS $U$ with an RRS $U'=\{u_1\dots u_{l-2},v,u_{l+1}\dots u_k\}$. 

 This first optimality condition implies that $\hat{w}_l$ cannot have a critical suffix $\xi$, as by \cite[Corollary 2.2]{Holt2012} we would then have $l[\tau(\xi)]=l[\tau(\hat{u}_l)]$ and we could have started the RRS at the first letter of $\xi$. This means that, assuming without loss of generality that $f[u_l]$ is positive, we must have  $p(\hat{w}_l)<p(\hat{u}_l)$.  Because $\beta_{l-1}=\epsilon$, we know that $u_l=f[u_l]w_l$, and so $n(\hat{w}_l)=n(\hat{u}_l)$ (Note here it is possible for $n(\hat{u})$ to be 0).

We also know that $\hat{v}$ is geodesic, so $p(\hat{v})+n(\hat{v})=m$. Together with the fact the $\hat{u}_{l-1}$ is a critical subword of $\hat{v}$ and $\hat{w}_{l}$ is a subword of $\hat{v}$, we get that $n(\hat{v})=n(\hat{u}_{l-1})=n(\hat{w}_l)$. Thus if $\hat{u}_{l-1}$ is positive then so are $\hat{w}_l$ and $\hat{u}_l$ and hence $\hat{v}$ will also be a positive critical word, and $l[\tau(\hat{v})]=l[\tau(\hat{u}_l)]$. Similarly, if $\hat{u}_{l-1}$ is unsigned, then $\hat{v}$ must start with an alternating positive sequence of length $p(\hat{v})$, and end in a negative alternating sequence of length $n(\hat{u}_l)=n(\hat{v})$. This means that $\hat{v}$ is a critical word of the same form as $\hat{u}_l$ and $l[\tau(\hat{v})]=l[\tau(\hat{u}_l)]. $

Repeating this operation for all problematic~$l$ will eventually produce an RRS satisfying all three optimality conditions.

\smallskip

Now we prove uniqueness.
 Suppose that $w=\mu w_1\cdots w_{k+1}\gamma$ and $w=\mu' y_1\cdots y_{k'+1}\gamma'$ are two factorizations of $w$ associated to two different optimal RRS $U=u_1,\dots, u_{k+1}$ and $V=v_1,\dots, v_{k'+1}$. Firstly, $\mu=\mu'$ because otherwise one of the RRS would not be optimal. 
 
 If $k=0$, this means that we can write $w=\mu a w'\gamma$, where~$a$ is a letter that commutes with all the letters in the word~$w'$ and $f[\gamma]=a^{-1}$. If $v_1$ were a P2G critical word beginning at $a$ then $v_1$ must contain a letter that does not commute with $a$, so $v_1$, considered as a subword of $w$, must end to the right of $f[\gamma]$. However, if $v_1$ were P2G critical, then $\hat{v_1}$ must be geodesic and so cannot contain $f[\gamma]$. This implies that~$v_1$ cannot be a P2G critical word.  Hence $k,k'=0$ with $u_1=w_1=y_1=v_1$.
 
 Assume then that $k,k'>0$. Observe that $f[w_1]=f[y_1]$ is a pseudo-generator of both $u_1$ and $v_1$. The other pseudo-generator of~$u_1$ and~$v_1$ must also coincide because this pseudo-generator is the closest letter to the right of~$f[w_1]$ that does not commute with~$f[w_1]$. So we need to prove $w_1=y_1$ and this will imply that $\tau(\hat{u}_1)=\tau(\hat{v}_1)$ and rest of the argument will follow by induction. 

Assume that $u_1\neq v_1$ and suppose without loss of generality that $v_1$ is a proper subword of $u_1$, so $u_1$ has nonempty intersection with $y_2$. Denote by $z$, the letter $l[\tau(\hat{v}_1)]$. The letter $z$ must be a pseudo-generator of both $u_1$ and $v_1$. 

\smallskip

First, if $k'=1$, then by \Cref{lemma:beta}, we know that $v_2=zy_2$, where all the letters in~$y_2$ commute with~$z$, but by optimality are not equal to $z$, and $f[\gamma']=z^{-1}$. However, $\hat{u}_1$ is geodesic word. This implies that $f[\gamma']$ is not in~$u_1$ and so~$l[u_1]$ must be in~$y_2$. Yet~$l[u_1]$ must either equal $z^\pm$ or be a pseudo-generator that does not commute with $z$. This means that we can assume that $k'>1$, and~$v_2$ must be a P2G critical word. 

Again, let  $z=l[\tau(\hat{v}_1)]^{\pm 1}$ be the pseudo-generator shared by~$v_1$ and~$v_2$. Let~$t_1$ be the other generator of~$\hat{v}_1$ and~$t_2$ be the other generator of~$\hat{v}_2$ (which may be equal to~$t_1$). The pseudo-generators of~$u_1$ are~$z$ and~$t_1$, so we now divide into cases according to whether $l[u_1]=z^\pm$ or $l[u_1]=t_1^{\pm 1}$. In either case, we will show that there are letters trapped in~$u_1$ or~$v_2$ contradicting the fact that these are P2G words. 

\smallskip 

Suppose that $l[u_1]=z^\pm$. The word $v_1$ is a subword of $u_1$, and in particular all the letters in $\beta_{v_1}$ appear in $u_1$  to the right of at least one instance of each of the pseudo-generators on $v_1$ and so must commute with $z$. By \Cref{lemma:beta}, $\beta_{v_1}\in\{t_2^\#,\epsilon\}$, and so the fact that $z$ does not commute with $t_2$ implies that $\beta_{v_1}=\epsilon$. By \Cref{remark:optimality}, the first instance of a pseudo-generator in $y_2$ must be $t_2^\pm$, so there must be an instance of $t_2^{\pm}$  to the left of $l[u_1]=z^\pm$ in the intersection of $u_1$ and $y_2$. Also $z^{\pm}$ appears in $v_1$. If $t_1\neq t_2$, then $t_2$ is an internal letter in $u_1$, but then it would be trapped in $u_1$ between instances of $z$ giving a contradiction. If $t_1=t_2$, by optimality $\alpha_{v_2}$ contains at least one letter that does not commute with both $t_1$ and $z$. Also the letters in $\alpha_{v_2}$ must appear before the first appearance of $t_2$ in $y_2$, so these letters are trapped in $u_1$ between either instances of $t_1$'s or instances of $z$'s.

\medskip

Suppose $l[u_1]=t_1^\pm$. Suppose that $l[u_1]$ does not lie in $y_2$. If $t_1\neq t_2$, then there is an instance of $t_2$ trapped in $u_1$ between instances of $z$'s. If $t_1=t_2$,  then, by optimality, $\alpha_{v_2}$ has at least one letter that either does not commute with $z$ or does not commute with $t_1$. This letter must be trapped in $u_1$.
Suppose that $l[u_1]$ lies in $y_2$. Then there is an instance of $t_1$ in $y_2$. Suppose $t_1\neq t_2$. If $z$ appears in $v_2$ to the right of $t_1$ then $t_1$ is trapped in $v_2$. Otherwise, $z$ must appear to the left of $t_1$ in the intersection of $u_1$ and $y_2$. Now, either we have that $\beta_{v_1}=t_2^\#$ or, by  Remark~\ref{remark:optimality}, the first instance of a pseudo-generator in $y_2$ must be~$t_2$, so in either case we have an instance of $t_2$ trapped in $u_1$ between instances of $z$. If $t_1=t_2$, by optimality, there exist letters in $\alpha_{v_2}$ that commute with~$z$ but not with $t_1$. Also, the letters of $\alpha_{v_2}$ appear in~$y_2$ before the first instance of~$t_1$, so they appear in $u_1$ before $l[u_1]=t_1$. Since there is an instance of $t_1$ in $v_1$, there are letters in $\alpha_{v_2}$ that are trapped in $u_1$. 
\end{proof}

%% file: RRS_kP2G.pdf_tex
\begingroup%
  \makeatletter%
  \providecommand\color[2][]{%
    \errmessage{(Inkscape) Color is used for the text in Inkscape, but the package 'color.sty' is not loaded}%
    \renewcommand\color[2][]{}%
  }%
  \providecommand\transparent[1]{%
    \errmessage{(Inkscape) Transparency is used (non-zero) for the text in Inkscape, but the package 'transparent.sty' is not loaded}%
    \renewcommand\transparent[1]{}%
  }%
  \providecommand\rotatebox[2]{#2}%
  \newcommand*\fsize{\dimexpr\f@size pt\relax}%
  \newcommand*\lineheight[1]{\fontsize{\fsize}{#1\fsize}\selectfont}%
  \ifx\svgwidth\undefined%
    \setlength{\unitlength}{408.12478854bp}%
    \ifx\svgscale\undefined%
      \relax%
    \else%
      \setlength{\unitlength}{\unitlength * \real{\svgscale}}%
    \fi%
  \else%
    \setlength{\unitlength}{\svgwidth}%
  \fi%
  \global\let\svgwidth\undefined%
  \global\let\svgscale\undefined%
  \makeatother%
  \begin{picture}(1,0.50116849)%
    \lineheight{1}%
    \setlength\tabcolsep{0pt}%
    \put(0,0){\includegraphics[width=\unitlength,page=1]{RRS_kP2G.pdf}}%
    \put(0.05495531,0.49420092){\makebox(0,0)[lt]{\lineheight{1.25}\smash{\begin{tabular}[t]{l}$\mu$\end{tabular}}}}%
    \put(0.32069323,0.49420093){\makebox(0,0)[lt]{\lineheight{1.25}\smash{\begin{tabular}[t]{l}$w_1$\end{tabular}}}}%
    \put(0.6706719,0.49420093){\makebox(0,0)[lt]{\lineheight{1.25}\smash{\begin{tabular}[t]{l}$w_2$\end{tabular}}}}%
    \put(0.91742726,0.49420093){\makebox(0,0)[lt]{\lineheight{1.25}\smash{\begin{tabular}[t]{l}$\gamma$\end{tabular}}}}%
    \put(0.17643561,0.4019436){\makebox(0,0)[lt]{\lineheight{1.25}\smash{\begin{tabular}[t]{l}$\alpha_1\rho_1$\end{tabular}}}}%
    \put(0.30924341,0.40194341){\makebox(0,0)[lt]{\lineheight{1.25}\smash{\begin{tabular}[t]{l}$\hat{u}_1$\end{tabular}}}}%
    \put(0.44319821,0.40194341){\makebox(0,0)[lt]{\lineheight{1.25}\smash{\begin{tabular}[t]{l}$\beta_1$\end{tabular}}}}%
    \put(0.30008494,0.31961696){\makebox(0,0)[lt]{\lineheight{1.25}\smash{\begin{tabular}[t]{l}$\tau(\hat{u}_1)$\end{tabular}}}}%
    \put(0.32286247,0.26836305){\makebox(0,0)[lt]{\lineheight{1.25}\smash{\begin{tabular}[t]{l}$l[\tau(\hat{u}_1)]$\end{tabular}}}}%
    \put(0.57642235,0.2193515){\makebox(0,0)[lt]{\lineheight{1.25}\smash{\begin{tabular}[t]{l}$u_2$\end{tabular}}}}%
    \put(0.38516611,0.15496359){\makebox(0,0)[lt]{\lineheight{1.25}\smash{\begin{tabular}[t]{l}$\alpha_2\rho_2$\end{tabular}}}}%
    \put(0.66693454,0.15496359){\makebox(0,0)[lt]{\lineheight{1.25}\smash{\begin{tabular}[t]{l}$\hat{u}_2$\end{tabular}}}}%
    \put(0.65548675,0.00216281){\makebox(0,0)[lt]{\lineheight{1.25}\smash{\begin{tabular}[t]{l}$\tau(\hat{u}_2)$\end{tabular}}}}%
    \put(0.76949429,0.05854209){\makebox(0,0)[lt]{\lineheight{1.25}\smash{\begin{tabular}[t]{l}$f[\gamma]^{-1}$\end{tabular}}}}%
    \put(0.89519339,0.00344419){\makebox(0,0)[lt]{\lineheight{1.25}\smash{\begin{tabular}[t]{l}$\gamma$\end{tabular}}}}%
    \put(0.87651244,0.05982357){\makebox(0,0)[lt]{\lineheight{1.25}\smash{\begin{tabular}[t]{l}$f[\gamma]$\end{tabular}}}}%
    \put(0,0){\includegraphics[width=\unitlength,page=2]{RRS_kP2G.pdf}}%
  \end{picture}%
\endgroup%

%% file: RRS_knotP2G.pdf_tex
\begingroup%
  \makeatletter%
  \providecommand\color[2][]{%
    \errmessage{(Inkscape) Color is used for the text in Inkscape, but the package 'color.sty' is not loaded}%
    \renewcommand\color[2][]{}%
  }%
  \providecommand\transparent[1]{%
    \errmessage{(Inkscape) Transparency is used (non-zero) for the text in Inkscape, but the package 'transparent.sty' is not loaded}%
    \renewcommand\transparent[1]{}%
  }%
  \providecommand\rotatebox[2]{#2}%
  \newcommand*\fsize{\dimexpr\f@size pt\relax}%
  \newcommand*\lineheight[1]{\fontsize{\fsize}{#1\fsize}\selectfont}%
  \ifx\svgwidth\undefined%
    \setlength{\unitlength}{408.12478854bp}%
    \ifx\svgscale\undefined%
      \relax%
    \else%
      \setlength{\unitlength}{\unitlength * \real{\svgscale}}%
    \fi%
  \else%
    \setlength{\unitlength}{\svgwidth}%
  \fi%
  \global\let\svgwidth\undefined%
  \global\let\svgscale\undefined%
  \makeatother%
  \begin{picture}(1,0.38937066)%
    \lineheight{1}%
    \setlength\tabcolsep{0pt}%
    \put(0,0){\includegraphics[width=\unitlength,page=1]{RRS_knotP2G.pdf}}%
    \put(0.03290321,0.38240311){\makebox(0,0)[lt]{\lineheight{1.25}\smash{\begin{tabular}[t]{l}$\mu$\end{tabular}}}}%
    \put(0.27658936,0.38240311){\makebox(0,0)[lt]{\lineheight{1.25}\smash{\begin{tabular}[t]{l}$w_1$\end{tabular}}}}%
    \put(0.6706719,0.38240311){\makebox(0,0)[lt]{\lineheight{1.25}\smash{\begin{tabular}[t]{l}$w_2$\end{tabular}}}}%
    \put(0.91742727,0.38240311){\makebox(0,0)[lt]{\lineheight{1.25}\smash{\begin{tabular}[t]{l}$\gamma$\end{tabular}}}}%
    \put(0.12865606,0.29014578){\makebox(0,0)[lt]{\lineheight{1.25}\smash{\begin{tabular}[t]{l}$\alpha_1\rho_1$\end{tabular}}}}%
    \put(0.34232131,0.29014559){\makebox(0,0)[lt]{\lineheight{1.25}\smash{\begin{tabular}[t]{l}$\hat{u}_1$\end{tabular}}}}%
    \put(0.27435769,0.20781914){\makebox(0,0)[lt]{\lineheight{1.25}\smash{\begin{tabular}[t]{l}$\tau(\hat{u}_1)$\end{tabular}}}}%
    \put(0.36814066,0.15789012){\makebox(0,0)[lt]{\lineheight{1.25}\smash{\begin{tabular}[t]{l}$l[\tau(\hat{u}_1)]=a$\end{tabular}}}}%
    \put(0.58180214,0.11255393){\makebox(0,0)[lt]{\lineheight{1.25}\smash{\begin{tabular}[t]{l}$u_2$\end{tabular}}}}%
    \put(0.89519339,0.00190757){\makebox(0,0)[lt]{\lineheight{1.25}\smash{\begin{tabular}[t]{l}$\gamma$\end{tabular}}}}%
    \put(0.87591863,0.05562411){\makebox(0,0)[lt]{\lineheight{1.25}\smash{\begin{tabular}[t]{l}$a^{-1}$\end{tabular}}}}%
    \put(0,0){\includegraphics[width=\unitlength,page=2]{RRS_knotP2G.pdf}}%
    \put(0.7001874,0.05404813){\makebox(0,0)[lt]{\lineheight{1.25}\smash{\begin{tabular}[t]{l}$l[\tau(\hat{u}_1)]=a$\end{tabular}}}}%
    \put(0,0){\includegraphics[width=\unitlength,page=3]{RRS_knotP2G.pdf}}%
    \put(0.65352817,0.15788862){\makebox(0,0)[lt]{\lineheight{1.25}\smash{\begin{tabular}[t]{l}$w_2$\end{tabular}}}}%
    \put(0.50066369,0.05524443){\makebox(0,0)[lt]{\lineheight{1.25}\smash{\begin{tabular}[t]{l}$w_2$\end{tabular}}}}%
  \end{picture}%
\endgroup%

%% file: Results.tex
We claim that in any $3$-free Artin group, any non-geodesic word can be reduced to an equivalent geodesic word, using only rightward reducing sequences and free reductions.

\begin{definition}
	We define the set $W$ to be the set of freely reduced words that admit no RRS.
\end{definition}

Notice that given a word~$w$, we may determine if $w$ is in the set~$W$ using a possibly very large, but still finite number of steps.  Because there are a finite number of subwords in $w$, there are a finite number of possible P2G critical subwords, and we must test whether each of these subwords might the first word in an RRS. {Nevertheless, we will give later a very efficient process to compute RRSs for our purposes}.

\medskip

We can now state a more detailed version of our main theorem:

\begin{maintheorem}\label{MT1}
	Let $A_S$ be an  3-free Artin group. Then $W$, the set of freely reduced words which do not admit an RRS, is the set of geodesic words of $A_S$. Moreover, membership in $W$ for a word $w$ can be determined with a finite time algorithm {that is quadratic on the length of $w$}, so the word problem is solvable for 3-free Artin groups {in quadratic time}. 
\end{maintheorem}

We will now outline the main steps of this theorem, with the proofs of propositions and auxiliary lemmas left to \Cref{Section:Propositions,Section:lemmas}. This strategy is similar to the one given for large Artin groups in \cite{Holt2012}. It is clear that words which admit an RRS are not geodesic, so the bulk of the proof of the main theorem involves showing that any non-geodesic word must admit an RRS.

\medskip

Recall that $\mathcal{A}^*$ is the set of all words in $S\cup S^{-1}$ where~$A_S$ is the Artin group generated by~$S$. 
  To show that any non-geodesic word admits an RRS, we would like to first examine non-geodesic words~$wg$ built out of a word $w\in W$ multiplied by a single letter $g\in S\cup S^{-1}$.
  
  \bigskip
 
\noindent \textbf{\Cref{P1}}. Let $A_S$ be a 3-free Artin group with generating set $S$.
	Let $w\in W$ and $g\in S\cup S^{-1}$ such that~$wg\not\in W$. Then there is a unique optimal RRS that can be applied to $wg$. Furthermore, this unique RRS applied to~$wg$ followed by a free reduction yields an element of~$W$.
	
\bigskip

{We also provide a linear-time algorithm to compute the optimal RRS of the previous proposition, which is much better than just looking at all the possible RRS in a naive way. This algorithm will be detailed in Section~5 and we provide explicit SAGEMath code for it in \cite{code}. Actually, in this program the user can enter a word in a 3-free Artin group and obtain a geodesic representative.}

\bigskip

\noindent \textbf{\Cref{algorithm}}. { Let $w\in W$ and $g\in S \cup S^{-1}$. The algorithm given below takes $w$ and $g$, decides if $wg$ admits an RRS and returns the optimal RRS $U$ and the corresponding decomposition of $wg$ into $\mu w_1w_2\dots w_{k+1} g$ in case it exists. Moreover, the algorithm complexity is linear on the length of $w$.}

{See the proof of \Cref{algorithm} for an explicit description of the algorithm.}
\bigskip

\Cref{P1,algorithm} imply that for any given word~$w$, we can find an equivalent word $\hat{\phi}(w)\in W$ which can be obtained from~$w$ via successive applications of free-reduction and RRS. We define a map $\hat{\phi}: \mathcal{A}^{*}\longrightarrow W$ {---which is inspired in the map $\rho$ of \cite{Holt2013}---} recursively as follows:
\begin{enumerate}
	\item For $w\in W$, $\hat{\phi}(w)=w$.
	\item For $w\in W, g\in S\cup S^{-1}$ and $wg\not\in W$, if $l[w]=g^{-1}$ then $\hat{\phi}(wg)=\mathrm{pref}[w]$ (i.e. the prefix of $w$ without its last letter), otherwise $\hat{\phi}(wg)$ is the element of $W$ obtained from $wg$ via an optimal RRS as specified in \Cref{P1}.
	\item For any $w\in \mathcal{A}^{*}$ and $g\in S\cup S^{-1}$, $\hat{\phi}(wg)=\hat{\phi}(\hat{\phi}(w)g)$.

\end{enumerate}

	We consider $W/\sim$ the quotient set of~$W$ under the following equivalence relation: $a,b\in W$ $a\sim b$ if $a,b$ have the same length and represent the same element in $G$. We consider the natural projection $\pi: W\longrightarrow W/\sim$ and define $\phi: \mathcal{A}^{*}\longrightarrow W/\sim$ such that $\phi= \pi\circ\hat{\phi}$.
	
In order to prove the main theorem, we would like to apply~$\phi$ to an arbitrary word~$w$ and compare~$\phi(w)$ to~$\phi(u)$ where $u$ is a geodesic representative of~$w$. If we can show that $\phi(w)=\phi(u)$, then $\hat{\phi}(w)$ will be a geodesic representative of the given group element. To that end, we would like to determine how the map $\phi$ behaves with respect to the insertion and deletion of subwords~$gg^{-1}$ for $g\in \mathcal{A}^*$ and how~$\phi$ behaves when changing a word using the Artin group relations.
	
\bigskip	
	
\noindent \textbf{\Cref{main_proposition}}. 
{Let $w$ be a word and be $w'$ the word obtained from $w$ by performing a relation $(a,b)_{m_{ab}}=(b,a)_{m_{ab}}$, inserting or deleting a subword $gg^{-1}$, $g\in S \cup S^{-1}$. Then $\phi(w)=\phi(w')$.} 

\bigskip

\begin{proof}[Proof of the main theorem]

{
First we note that an RRS followed by a free reduction shortens a word and therefore any word not in $W$ cannot by geodesic. }

{To show the converse, consider some $w\in W$, so $\hat{\phi}(w)=w$. Suppose that $u$ is a geodesic representation of the same group element. 
We can define a chain of words $w=w_0$, $w_1$,...,$w_k=u$, where for each $i=1,..,k$, $w_i$ is obtained from $w_{i-1}$ either by insertion or deletion of a subword $gg^{-1}$ for $g\in S\cup S^{-1}$ or by the replacement of a subword $(a,b)_{m_{ab}}$ by $(b,a)_{m_{ab}}$. }

{Using \Cref{main_proposition}, we get that  $\phi(w)=\phi(w_1)=\phi(w_2)=\dots \phi(w_k)=\phi(u)$. 
Thus $\hat{\phi}(w)$ and $\hat{\phi}(u)$ have the same length. Finally $u$ is a geodesic word and therefore $u\in W$ and $\hat{\phi}(u)=u$. Thus $u$ and $w$ have the same length and so $w$ is also geodesic.}

{In \Cref{algorithm}, we give a specific algorithm to decide the existence and compute the optimal RRS at every step of $\hat{\phi}$, and we prove that this algorithm has linear complexity on the length of the input word. Since we need to run this algorithm at most once for each letter in $w$, we have a quadratic time algorithm to solve the word problem.}

\end{proof}

The main theorem shows that any word can be reduced to a geodesic representative using only RRS and free reductions. However geodesics are not unique in Artin groups, so we would also like to know how different geodesics representing the same element might be related via RRS. 

\begin{corollary}\label{Cor:LengthPreserved}
	Let $v_1,v_2$ be two geodesic words representing the same element in a 3-free Artin group $A$. Using only commutations and $\tau$-moves, we can transform~$v_1$ into~$v_2$ without increasing the length at any step.
\end{corollary}

\begin{proof}
	We will prove this by induction on the length of~$v_1$ (and $v_2$). For length 1 the result follows. Suppose that the length of~$v_1$ and~$v_2$ is~$n$. Our induction hypothesis states that for any pair of geodesics $u_1=_G u_2$ of length less than~$n$, we can transform~$u_1$ into~$u_2$ without increasing the length. 
	
	Let $a=l[v_1]$, then $v_1=w_1 a$ where~$w_1$ is geodesic of length~$n-1$. Note that $v_2a^{-1}$ is not geodesic because~$w_1$ is a shorter representative of the same element. By \Cref{P1} there is an RRS from $v_2a^{-1}$ to a word $w_2=_G w_1$. This means that there is an RRS within $v_2$ cancelling~$a$, so we can transform~$v_2$ into~$w_2a$ using only commutations and $\tau$-moves. This proposition also implies that~$w_2$ is geodesic of length~$n-1$, so applying the induction hypothesis we can change~$w_1a$ into~$w_2a$ without increasing the length.
\end{proof}

\subsection{Property H} \label{Sec:PropH}

In \cite{Dehornoy2011}, Dehornoy defined a group  with finite positive presentation $\langle S|R \rangle$ that satisfies \textit{Property H} as a group where any word representing the identity element, $w=_G\epsilon$, can be reduced to the empty word via only 3~types of operations. 

\begin{description}
\item[Type 0] Free reductions.
\item[Type 1] Replacing $v$ with $u$ or $v^{-1}$ with $u^{-1}$ where $v=u$ is a relation in $R$. 
\item[Type 2] If $u_1u_2$ and $v_1v_2$ are two words representing the same group element then we can replace $v_1^{-1}u_1$ with $v_2u_2^{-1}$ (here we require $u_1,v_1\neq \epsilon$) or $u_2v_2^{-1}$ with $u_1^{-1} v_1$ (with $u_2,v_2\neq \epsilon)$.
\end{description}

Property H can be thought of as a weak hyperbolicity property and is conjectured to be true for all Artin groups. 

\begin{corollary} [Property~H]
Let $A_S$ be a 3-free Artin group. Then this group satisfies Property~H.
\end{corollary}

\begin{proof}
    Godelle and Rees \cite{GodelleRees2016} show that $\tau$-moves in dihedral Artin groups are composed of repeated applications of operations of type~1 and~2. (See Lemmas 2.11 and 2.12 of that paper). Because our algorithm for reducing a word the identity involves only  $\tau$-moves applied to dihedral subwords, commutations (type~1 moves), and free reductions, this algorithm shows that 3-free Artin groups satisfy Property~H. 
\end{proof}

\subsection{Why 3-free?}

Here, we give further details about why the 3-free condition is stipulated in our algorithm. First, if relations of length 3 are allowed, we can find examples of non-geodesic words that cannot be reduced to geodesic words using this algorithm. 

\begin{example} Consider a braid group on 4 strands with generators $a,b,c$ where $m_{a,b}=m_{b,c}=3$ and $m_{a,c}=2$. The word $cbabc^{-1}ba^{-1}$ is not a geodesic word, as can be seen by repeated applications of $\tau$ moves to P2G critical words. 

\[c(bab)c^{-1}ba^{-1}\overset{\tau}{\to}(cabac^{-1})ba^{-1}\overset{\tau}{\to}ab^{-1}c(bab)a^{-1}\overset{\tau}{\to}ab^{-1}cabaa^{-1}\]

However when we put these $\tau$ moves together they do not form an RRS. The problem is that when we allow for relations of length 3, we can have P2G critical words nested inside each other. This means that there are many more possible sequences in which to apply $\tau$ moves to a given word. Without a canonical way to choose one sequence over another, it becomes more difficult to show that every non-geodesic word admits a finite such sequence. \end{example}

\medskip

The 3-free condition also gives easily verifiable  conditions on a word that imply there are no P2G critical subwords and thus no $\tau$-moves are possible. 

\begin{example}
Suppose that $a,b,c$ are letters such that $c$ commutes with $a$, but not $b$. If $m_{a,b}=4$ then $abcab$ cannot be P2G critical because $c$ is trapped in between copies of $b$. However if $m_{a,b}=3$, then $acbca$ is a P2G critcal word. This comes from that fact that as there is only one copy of the letter $b$ needed to make an alternation of length 3, and thus it is much harder to trap letters. 
\end{example}

We use the argument that a letter is trapped frequently when showing that an RRS with a particular factorization is not possible, so it is hugely advantageous to have weak conditions which imply a letter is trapped.

%% file: Props.tex
The aim of this section is to prove the propositions stated earlier and required for the proof of the main theorem. Auxiliary results that will be proven in \Cref{Section:lemmas} are now stated.

\bigskip

%



\noindent \textbf{\Cref{Lem:commutation does not affect RRS}}   Suppose that $A_S$ is a 3-free Artin group. 
Consider a word $w$ admitting an RRS $U= u_1,\dots, u_{k+1}$ with factorization $w=\mu w_1\dots w_{k+1} \gamma$. Suppose also that  $w$ has $ab$ as a subword, where $a,b$ are letters such that $m_{a,b}=2$. Let $w'$ be the word obtained from $w$ by replacing the subword $ab$ with the subword $ba$. 

\medskip\noindent
a) Then $w'$ also admits an   RRS $V= v_1,\dots, v_{n+1}$ with factorization $w'=\nu y_1\dots y_{n+1}\eta$, where $f[\gamma]=f[\eta]$, $n=k$, and $\hat{u}_i=\hat{v}_i$ for $i\leq k$. 

\medskip\noindent
{b) Moreover, if $ab$ does not intersect with $\gamma$, $U$ takes $w$ to $u \mathrm{suf}[\gamma]$ and $V$ takes $w'$ to $v \mathrm{suf}[\gamma]$, then we can transform $v$ into $u$ by applying only commutations.} 

\medskip\noindent
{c) In the hypothesis of b), if $\mathrm{pref}[w]\in W$ and $U$ is optimal, then $V$ is also optimal. }

\bigskip

\noindent \textbf{\Cref{Lem: tau moves do not affect RRS}}  Suppose that $A_S$ is a 3-free Artin group. 
Consider a word~$w$ admitting an optimal RRS $U= u_1,\dots, u_k$ with factorization $w=\mu w_1\dots w_{k+1} \gamma$. Suppose also that  $w$ has a 2-generated critical subword $\zeta$, in the generators $s$ and $t$, such that $m_{s,t}>3$. Let~$w'$ be the word obtained from~$w$ by replacing the subword $\zeta$ with the subword $\tau(\zeta)$.

\medskip\noindent
a) $w'$ also admits an RRS $V=v_1,\dots, v_{n+1}$ with factorization $w'=\nu y_1\dots y_{n+1}\eta$.

\medskip\noindent
{b) Also, if ~$\zeta$ does not intersect with~$\gamma$, then $\gamma=\eta$. In this case, if $U$ takes $w$ to $u \mathrm{suf}[\gamma]$ and $V$ takes $w'$ to $v \mathrm{suf}[\gamma]$, then we can transform $v$ into $u$ by applying only commutations and $\tau$-moves.}

\medskip\noindent
{ c) In the hypothesis of b), if $\mathrm{pref}[w]\in W$, then $V$ is also optimal. }

\bigskip

\begin{proposition}\label{P1} Let $A_S$ be a 3-free Artin group with generating set $S$.
	Let $w\in W$ and $g\in S\cup S^{^{-1}}$ such that $wg\not\in W$. Then there is a unique optimal RRS that can be applied to $wg$. Furthermore, this unique RRS applied to $wg$ followed by a free reduction yields an element of $W$.
\end{proposition}

\begin{proof} 
By \Cref{ExistenceOfOptimalRRS}, $wg$ admits a unique optimal RRS. We apply the sequence of $\tau$-moves in this RRS to obtain an equivalent word of the form $w'g^{-1}g$. We would like to show that $w'$ does not admit an RRS, so we suppose otherwise. If $w'$ admits an RRS, then so does $w'g^{-1}$ (with a slightly longer tail). We can use \Cref{Lem:commutation does not affect RRS,Lem: tau moves do not affect RRS} to move $g^{-1}$ back into $w'$, changing the word $w'g^{-1}$ back into the word $w$ and maintaining the existence of an RRS at each step. Thus we contradict the original assumption that $w$ did not admit an RRS. 
\end{proof}

\begin{proposition}\label{algorithm} Let $w\in W$ and $g\in S \cup S^{-1}$. {We describe below an explicit} algorithm that takes $w$ and $g$, decides if $wg$ admits an RRS and returns the optimal RRS $U$ and the corresponding decomposition of $wg$ into $\mu w_1w_2\dots w_{k+1} g$ in case it exists. The complexity of this algorithm is linear in the length of $w$.
\end{proposition}

\begin{proof} 

 We will explain step by step how to factor $w$ into $w=\mu w_1\dots w_{k+1}$. We will construct every $w_i$ by starting from the last letter of $w$ and going from right to left. {We already know that if $wg$ admits an RRS, it admits a unique optimal RRS by \Cref{ExistenceOfOptimalRRS}. Optimality is strong condition and by \Cref{ExistenceOfOptimalRRS} if we cannot satisfy these conditions, we may conclude that $wg$ does not admit an RRS.}

Recall that for a given word $w_i$ and generators $a,c$ we can compute  $m(\hat{w}_i)$ by deleting all the letters not equal to $a^\pm$ or $c^\pm$, to obtain $\hat{w}_i$ in the dihedral Artin group generated by $a$ and $c$ and counting $m(\hat{w}_i)=n(\hat{w}_i)+p(\hat{w}_i)$. {Because we do not only want to check this number but also the form required by the criticality of the word, we will define $m'(\hat{w}_i)$ as follows: Suppose without lost of generality that the last letter of $\hat{w}_i$ is positive (the other case is analogous). If $\hat{w}_i$ is positive then $m'(\hat{w}_i)=\max\{r\mid {}_{r}(a,c) \text{ is a prefix or suffix of } \hat{w}_i \}$. Otherwise we set $\{x,y\}:=\{a,c\}$ and then $m'(\hat{w}_i)=\max\{r+r'\mid {}_{r}(a^{-1},c^{-1}) \text{ is a prefix  and } _{r'}(x,y) \text{ is a suffix of } \hat{w}_i \}$.}

Let $h$ be the rightmost letter of $w$. 

\underline{To build $w_{k+1}$:}
\begin{enumerate}
\item While $h$ commutes with $g$ but is not equal to $g^{-1}$, add $h$ to the left of $w_{k+1}$, and proceed to the next letter on the right on $w$. 
\item If $h=g^{-1}$, set $k=0$, and let $h$ be the first letter of $w_1$. The algorithm finishes.
\item {If $h$ commutes with all the letters in $w$ and $h$ is never equal to $g^{-1}$, there is no RRS for $wg$.}
\item Otherwise $h$ is a letter that does not commute with $g$ (Note $h$ cannot equal $g$ as this would violate the assumption that $w\in W$). Add $h$ to $w_k$. This letter cannot belong to $w_{k+1}$, so it is $l[w_k]$ and the pseudo-generators of $u_k$ are $h$ and $g$. Set $a$ to be the name of the letter $h$, $c$ to be the name of the letter $g$,  $m:=m_{a,c}$ and $i:=k$.
\end{enumerate}


\medskip
\underline{To build $w_i$}
 Let $h$ be the rightmost remaining letter of $w$. {If at some point $w$ is empty and the algorithm has not finished, then there is no RRS.}
\noindent
\begin{enumerate}
\item \label{start_w_i} While $m'(\hat{w}_i)<m-2$, add $h$ as a prefix to $w_i$ and set $h$ to be the next rightmost letter in $w$, and repeat this step. Note by \Cref{lemma:beta} and 3-freeness, $\beta_{i-1}$ is a power of either $a$ or $c$, so to create $u_i$ from $w_i$ we add a prefix which increases {$m'(\hat{w}_i)$} by at most 2. This means that while {$m'(\hat{w}_i)<m-2$}, we must continue to add letters to $w_i$ in order to end up with a $u_i$ as a P2G critical word.  {If there is a letter added during this process it has to be equal to $a^\pm$, $c^\pm$ or be a letter in $\beta$ or to~$\rho$. If this is not the case, then there is no RRS. }

\item \label{step: m-2} While $m'(\hat{w}_i)=m-2$:
\begin{enumerate}
\item If $h\in\{a^{\pm},c^\pm\}$ or if $h$ commutes with both $a$ and $c$ add it as a prefix to $w_i$, set $h$ to be the next rightmost letter in $w$ and go back to step \ref{start_w_i}.  Because {$m'(\hat{w}_i)=m-2$}, in order to build a P2G critical word, we either need to make $w_i$ longer or we must have $\beta_{i-1}\neq \epsilon$. By \Cref{lemma:beta} and 3-freeness, any non-empty $\beta_{i-1}$ must be a power of $a$ or $c$, and must commute with but be unequal to $l[w_{i-1}]$. We conclude that $h\neq l[w_{i-1}]$ in this case, and we must add $h$ to $w_i$.
\item If $h=b\notin\{a^{\pm},c^\pm\}$ and commutes with only with one of the letters in $\{a,c\}$,  let $l[w_{i-1}]=h$. Let $c_1\in\{a,c\}$ be the letter that commutes with $b$ and $c_2\in\{a,c\}$ be the letter that does not. Then the pseudo-generators of $\hat{u}_{i-1}$ are $c_2$ and $b$, and  $\beta_{i-1}$ will be some power of $c_1$. We can continue from Step \ref{start_w_i} adding letters to $w_{i-1}$ by setting $i:=i-1$, $m=m_{c_2,b}$, $a:=b$ and $c:=c_2$, and letting $h$ be the next rightmost letter in $w$.
\item {Otherwise, we encounter a situation where $h$ does not commute with $a$ and does not commute with $c$ This implies that $h$ would be trapped, and there will not be a way to build a P2G critical word $u_i$, so there is no RRS.}
\end{enumerate}

\item\label{step:m-1} While $m'(\hat{w}_i)=m-1$: Let $C\subset\{a^\pm, c^\pm\}$ be the set of possible missing letters to make $\hat{u}_i$ critical. Notice that the missing letter does {not} need to be unique, however $C$ must be a subset of either $\{a^\pm\}$, or $\{c^\pm\}$. {If {the part of $w_i$ that we have already constructed} begins in $a$, then $a^{-1}$ cannot be the missing letter as then the word $\hat{w}_i$ would not be critical, and $a$ cannot be the missing letter, as this would violate the optimality of the RRS {(\Cref{remark:optimality})}. Thus the set $C$ contains at most two letters which must have the same name, {say $c_1$}.}

\begin{enumerate}
\item \label{step:done} If $h\in C$, add $h$ to $w_i$, set $i=1$ and the algorithm finishes. 
\item If $h\not\in C$ and $h\in\{a^\pm,c^\pm\}$, add $h$ to $w_i$, let $h$ be the next rightmost letter in $w$ and go back to Step~\ref{start_w_i}. At this point {the only letters we could possibly have added to $\alpha_i$ commute with both pseudo-generators of $u_i$}. Thus by the third optimality condition, we know that {either we must add more letters to $\alpha_i$} {that do not commute with one of the generators} or the pseudo-generators of $u_{i-1}$ are not $\{a,c\}$. In either case, $h$ cannot be $l[w_{i-1}]$. 
\item If $h$ commutes with both $a^\pm$ and $c^\pm$  but is not equal to any of then, then add it to $w_i$, set $h$ to be the next rightmost letter in $w$ and go back to Step~\ref{start_w_i}. Either $a$ or $c$ must be a pseudo-generator of $u_{i-1}$ so a letter commuting with them both cannot be $l[w_{i-1}]$. 

\item \label{step: i-1}If $h=b\not\in\{a^\pm,c^\pm\}$ and does not commute with $c_1$ then {$h$} must be $l[w_{i-1}]$.
Otherwise it would be trapped in $u_i$. Continue completing $w_{i-1}$ from Step~\ref{start_w_i} by doing $i:=i-1$, $m=m(b,c_1)$, $c:=c_1$, $a:=b$.  
 
\item \label{alpha} If $h$ commutes with some letter $c_1$ in $C$ but not with the other letter $c_2$ in $\{a,c\}$, add it to $w_i$. This means that $\alpha_i\neq \epsilon$ {and contains some letter not commuting with a generator}. Let $h$ be the next rightmost letter in $w$. Do not go back to Step~\ref{start_w_i} yet. 

\begin{enumerate}

\item If $h=c_1$, add $h$ to $w_i$, set $i=1$ and the algorithm finishes.

\item If $h\neq c_1^\pm$ and commutes with $c_1$, add it to $w_i$, as a letter in $\alpha_i$.  Let $h$ be the next rightmost letter in $w$ and go to step \ref{alpha}.

\item If $h\not\in\{a^\pm,c^\pm\}$  and does not commute with $c_1$, go to step \ref{step: i-1}, and follow the instructions to start building $w_{i-1}$.

\item If $h=c_2^\pm$, then, since $\alpha_i\neq \epsilon$, $h=l[w_{i-1}]$. We continue completing $w_{i-1}$ from Step~\ref{start_w_i} by setting $i:=i-1$. In this case the pseudo-generators of $u_i$ and $u_{i-1}$ are the same.  

\item {If $h=c_1^{-1}\not\in C$, we will need to add some instance of $c_2^\pm$ to make $\hat{w}_i$. But, since $\alpha$ is not trivial, this is not possible, so there is no RRS.}

\end{enumerate}
\end{enumerate}
 
\item\label{final_step} {Run the RRS. If at some point $\hat{u}_i$ is not critical, then, since our previous choices to construct the RRS were unique, there is no RRS.}
\end{enumerate}

{Notice the algorithm is just a letter by letter verification of a finite number of conditions. This number does not depend on any variable (checking the $m'$ conditions just need comparing the new letter with the previous one and having two counters for $r$ and $r'$). Also, the complexity of running an RRS is linear on the length of $w$. So the complexity of the algorithm is linear on the length of~$w$. For an explicit code for this algorithm in SAGEMath see \cite{code}.}
\end{proof}

\medskip

{In \Cref{P1}, we showed that $\hat{\phi}$ was a well-defined map, and in \Cref{algorithm} we described an algorithm that can be used to recursively compute $\hat{\phi}(w)$. Now, we can state further auxiliary results that will be used in the proofs of the subsequent propositions. Again, for proofs of these lemmas, see \Cref{Section:lemmas}.}

\bigskip

\noindent \textbf{\Cref{Lemma_Case_ab=ba}} 
 {
Suppose that $w\in W$, $m_{a,b}=2$ and in the process of computing $\hat\phi(wab)$ we apply two RRS, one cancelling $a$ and another cancelling $b$. Suppose $wb\not\in W$. Then $\hat{\phi}(wba)$ can be obtained from  $\hat{\phi}(wab)$ via only commutations and $\tau$-moves.   }

\bigskip

\noindent \textbf{\Cref{P3}}
	$\hat\phi(w(a,b)_{m_{ab}})=\hat\phi(w(b,a)_{m_{ab}})$ for any $w\in W$, $a,b\in S$ and where $m_{ab}>2$ is the length of the relation involving $a$ and $b$.

\bigskip

\begin{proposition}\label{main_proposition}
{Let $w$ be a word and be $w'$ the word obtained from $w$ by performing a relation $(a,b)_{m_{ab}}=(b,a)_{m_{ab}}$, inserting or deleting a subword $gg^{-1}$, $g\in S \cup S^{-1}$. Then $\phi(w)=\phi(w')$.} 
\end{proposition}

\begin{proof}
Suppose that we are in the case where~\( w' \) is obtained from~\( w \) by applying the relation \( (a,b)_{m_{ab}} = (b,a)_{m_{ab}} \). Notice that, since \( \hat{\phi} \) is defined recursively (letter by letter), the fact that \( \hat{\phi}(w(a,b)_{m_{ab}}) = \hat{\phi}(w(b,a)_{m_{ab}}) \) holds is sufficient to obtain our result. Therefore, the case \( m_{ab} > 2 \) is proven by \Cref{P3}.

Suppose $m_{a,b}=2$ and we have performed a commutation $ab=ba$ in $w=w_1abw_2$. Let $U$ be the first optimal RRS applied to compute $\hat{\phi}(w)$ that cannot be applied in the computing of $\hat{\phi}(w')$ (because it involves the letters of the commutation or some letter to the right of them).  

Assume that the letter erased by $U$ is after $ab$. By \Cref{Lem:commutation does not affect RRS}, in the computing of $\hat{\phi}(w')$ we can apply an optimal RRS $V$ that cancels the same letter as $U$.
Moreover, if $U$ produces the word $u$ and $V$ produces the word $v$, the lemma tells us that we can transform $u$ to $v$ by using commutations that do not affect the tails of the RRSs. This means that if we continue computing $\hat{\phi}(w)$ and we encounter an RRS, we can use \Cref{Lem:commutation does not affect RRS} and repeat the process in the same manner. In the end, we have reduced the same number of letters in each word so $\phi(w)=\phi(w')$. 

Assume that the letter erased by $U$ is $b$, then clearly $\hat\phi(\hat{\phi}(w_1)ab)= \hat\phi(\hat{\phi}(w_1)ba)$, because $a$ is in the last factor $w_{k+1}$ of $U$.
Suppose that letter erased by $U$ is $a$. If $\hat{\phi}(w_1)b\in W$, then the case is analogous to the previous one. Otherwise, \Cref{Lemma_Case_ab=ba} shows that $\hat{\phi}(wab)$ is equivalent to $\hat{\phi}(wba)$ by commutations and $\tau$-moves. We continue the process using \Cref{Lem:commutation does not affect RRS,Lem: tau moves do not affect RRS} as in the previous case and we also get $\phi(w)=\phi(w')$.

Suppose that we have inserted a subword $gg^{-1}$ in $w=w_1w_2$ obtaining $w'=w_1gg^{-1}w_2$. If $\hat{\phi}(w_1)g$ does not admit an RRS, the result follows. Otherwise, consider the optimal RRS in $\hat{\phi}(w_1)g$, that must cancel the last letter, and therefore transform $\hat{\phi}(w_1)$ into $w_1'g^{-1}$. {We claim that t}here is no RRS in $w'_1g^{-1}$ cancelling $g^{-1}$. For there {must be} a sequence of commutations and $\tau$-moves that reverse the original RRS taking $\hat{\phi}(w_1)g\overset{RRS}\rightsquigarrow w'_1g^{-1}g$. Thus by \Cref{Lem:commutation does not affect RRS,Lem: tau moves do not affect RRS}, we see that the existence of an RRS in $\hat{\phi}(w_1)g^{-1}$ implies the existence of an RRS in {$w_1$} contradicting {our definition of $\hat{\phi}$ wich is applies recursively left to right}. 

{This implies that $w'_1g^{-1}\in W$. To continue computing $\hat{\phi}(w')$, we must compute $\hat{\phi}(\hat{\phi}(w_1gg^{-1})w_2)$.} The next optimal RRS in the computing of $\hat{\phi}(\hat{\phi}(w_1)w_2)$ may be not applicable in $\hat{\phi}(\hat{\phi}(w_1gg^{-1})w_2)$. But again, since $\hat{\phi}(w_1)$  is equivalent to  $\hat{\phi}(w_1gg^{-1})$ by commutations and $\tau$-moves, we can also find an optimal RRS $V$ that cancels the same letter as $U$. Moreover, if $U$ produces the word $u$ and $V$ produces the word $v$, this lemma tells us that we can transform $u$ to $v$ by using commutations and $\tau$-moves that do not affect the tails of the RRSs. This allow us to reapply \Cref{Lem:commutation does not affect RRS,Lem: tau moves do not affect RRS} every time we encounter an RRS. In the end, we have reduced the same number of letters in each word so $\phi(w)=\phi(w')$.

\end{proof}

%% file: Lemmas.tex
\begin{lemma}\label{Lem:commutation does not affect RRS}  Suppose that $A_S$ is a 3-free Artin group. 
Consider a word $w$ admitting an RRS $U= u_1,\dots, u_{k+1}$ with factorization $w=\mu w_1\dots w_{k+1} \gamma$. Suppose also that  $w$ has $ab$ as a subword, where $a,b$ are letters such that $m_{a,b}=2$. Let $w'$ be the word obtained from $w$ by replacing the subword $ab$ with the subword $ba$. 

\medskip\noindent
a) Then $w'$ also admits an   RRS $V= v_1,\dots, v_{n+1}$ with factorization $w'=\nu y_1\dots y_{n+1}\eta$, where $f[\gamma]=f[\eta]$, $n=k$, and $\hat{u}_i=\hat{v}_i$ for $i\leq k$. 

\medskip\noindent
{b) Moreover, if $ab$ does not intersect with $\gamma$, $U$ takes $w$ to $u \mathrm{suf}[\gamma]$ and $V$ takes $w'$ to $v \mathrm{suf}[\gamma]$, then we can transform $v$ into $u$ by applying only commutations.} 

\medskip\noindent
{c) In the hypothesis of b), if $\mathrm{pref}[w]\in W$ and $U$ is optimal, then $V$ is also optimal. }

\end{lemma}

\begin{proof}

We divide into various cases according to where the subword $ab$ appears in the factorization $\mu w_1\dots w_k w_{k+1}\gamma$. {In each case, we first describe how to obtain a new RRS, $V$, from the given RRS, U. We then verify that the new RRS V satisfies the conditions required by b) above. Finally to prove c), in each case, we verify that the algorithm given in \Cref{algorithm} applied to the word $w'$ will give the RRS $V$, and thus $V$ must be optimal.}

\underline{If $ab$ is in $\mu$, or $ab$ appears to the right of $f[\gamma]$:} then this is immediate. {In the first case, this commutation takes $u$ to $v$ {and V is obviously optimal}. In the second $\eta\neq \gamma$.}

\underline{If $a=l[\mu]$ and $b=f[w_1]$:} The word $w'$ admits an RRS $v_1,\dots, v_{n+1}$, with $v_i=u_i$ if $i>1$ and where $v_1$ is constructed from $u_1$ by inserting the letter~$a$ immediately to the right of the letter $f[u_1]=b$. In the reduction process, the letter~$a$ will end up in $\alpha_{v_1}$, and so inserting this letter does not affect the RRS to the right of~$u_1$. {In this case $u=v$.}  {To check c), notice that the algorithm in \Cref{algorithm} works the same for $w'$ until it arrives to $a$, then we have to apply case 3.c) + 3.a) or case 3.e.i). The algorithm finishes and gives exactly $V$.}

\underline{If $a=f[w_1]$:} In a similar argument to the previous case, we see that the word $w'$ admits an RRS $v_1,\dots, v_{n+1}$, with $v_i=u_i$ if $i>1$ and where~$v_1$ is constructed from $u_1$ by deleting the letter~$b$ immediately to the right of $f[u_1]=a$. {In this case $u=v$.} {Notice that, when arriving to $a$, the algorithm in \Cref{algorithm} finishes using case 3.a) or 3.e.i) and gives exactly $V$.}

\underline{If $ab$ is contained in $w_i$ for $i\leq k$ {but is neither a prefix nor a suffix of $u_i$}:} At most one of the letters $a,b$ can be a pseudo-generator for this critical word~$u_i$. Let~$t$ be a letter in $\{a,b\}$ that is internal in~$u_i$. The letter~$t$ can pass via commutations to either the left or the right of the word~$\hat{u}_i$. If~$t$ can pass to the left then~$t$ ends up in~$\alpha_{u_i}$ or~$\rho_{u_i}$ and replacing~$ab$ with~$ba$ in~$u_i$, will not change this fact. Similarly, if~$t$ cannot pass via commutations to the left of $\hat{u}_i$ but can pass to the right and~$t$ ends up in $\beta_{u_i}$, then replacing~$ab$ by~$ba$, $t$ still cannot pass to the left but can pass to the right. Also note that by Lemma \ref{lemma:beta}, we know that $\beta_{u_i}$ cannot contain both~$a$ and~$b$, so at most one of these letters is passed to the right. 

Thus, starting with the word~$u_i$ and replacing~$ab$ with~$ba$, we obtain a P2G critical word~$v_i$, where $ \hat{u}_i=\hat{v}_i$. We also see that $\beta_{u_i}=\beta_{v_i}$. {We have $u=v$ except in the case where both~$a$ and~$b$ end up in either~$\alpha_i$ or~$\rho_i$. In such a case, applying again $ba=ab$ will take $v$ to $u$.} {To see c), we apply the algorithm in \Cref{algorithm}.  Since $b\neq l[u_i]$, the algorithm tells us that the pseudo-generators of~$u_{i}$ are the same as in the corresponding~$v_{j}$ ($j=i+(n-k)$). This plus $a\neq f[u_i]$ implies that $ba$ will be incorporated into~$y_j$. For the rest of the algorithm the process will be exactly the same as in~$U$, obtaining~$V$ (so $i=j$).}

\underline{If $ab$ is a suffix of $w_i$ where $i\leq k$}: In this case, $b$ is a pseudo generator of $u_i$, and $a$ is not. We obtain a new RRS with $v_i$ obtained from $u_i$ by deleting the last instance of $a$ in $u_i$, and letting $y_{i+1}=aw_{i+1}$.  In this new RRS, we have $\hat{v}_i=\hat{u}_i$. 

We also know that either $a$ can pass via commutations to the left of $\hat{u}_i$, {and $a= l[\rho_{u_i}]$}, or $a$ cannot pass to the left via commutations and {$\beta_{u_i}=a^\#$}. In the first case, $a$ commutes with all the letters of $l[\tau(\hat{u}_i)]\beta_{u_i}=l[\tau(\hat{v}_i)]\beta_{v_i}$ and in $w'$
 we have $v_{i+1}=l[\tau(v_i)]\beta_{v_i}aw_{i+1}$, so~$a$ can pass to the left via commutations ending up in either~$\alpha_{v_{i+1}}$ or~$\rho_{v_{i+1}}$, leaving $\hat{v}_{i+1}=\hat{u}_{i+1}$. {Also, we can use further commutations to push $a$ to the left, transforming $v$ into $u$}. In the second case we have that $\beta_{u_i}=\beta_{v_i}a$ and we end up with $v_{i+1}=u_{i+1}$ {and $u=v$}.
 
 {To check c), we apply the algorithm in \Cref{algorithm} to $w'$. If $a^\#=\beta_{u_i}$, then: arrive to $a$, apply step 2.a) and incorporate $a$ to $y_{j+1}$ ($j=i+(n-k)$);  arrive to $b$, apply 3.d) and incorporate $b$ to $y_j$; continue as the algorithm for $U$ and we obtain $V$ (so $i=j$). If $a= l[\rho_{u_i}]$, then: arrive to $a$ (which cannot be a generator of $v_{j+1}$) and $b$, apply step 2.a) + 2.b), or 3.c)+ 3.d); continue as in $U$ and obtain $V$.}

\underline{If $a=l[w_i]$ and $b=f[w_{i+1}]$ for $1\leq i<k$:} In this case we know that~$u_i$ is P2G critical and that~$a$ is a pseudo-generator of~$u_i$. We define a new RRS on~$w'$ by letting~$y_i$ be the word created from~$w_i$ by inserting the letter b just before the last letter of~$w_i$. We let~$y_{i+1}$ be the word obtained from~$w_{i+1}$ by deleting the first letter. Using this new factorization of~$w'$, we obtain an RRS. Because~$b$ commutes with $a=l[y_i]$, it cannot be a pseudo-generator of~$v_i$, we have that $\hat{v}_i=\hat{u}_i$.

In the original RRS $u_{i+1}=l[\tau(\hat{u}_i)]\beta_{u_i}w_{i+1}$ and~$b$ is the first letter of~$w_{i+1}$. If $b$ is internal in $u_{i+1}$, then it commutes with both generators of $u_i$ and this instance of~$b$ can be pushed via commutations to the left to become part of~$\rho_{u_{i+1}}$. The removal of $b$ from $w_{i+1}$ does not affect the remainder of the RRS and therefore $f[\eta]=f[\gamma]$. 
{After the RRS $U$, we end up with a word containing $ \mathrm{pref}[\tau(\hat{u}_i)]\alpha_{u_{i+1}}\rho_{u_{i+1}} \mathrm{pref}[\tau(\hat{u}_{i+1})].$ The instance of $b$ in the beginning of $w_{i+1}$ ends up in $\alpha_{u_{i+1}}\rho_{u_{i+1}}$. 
Removing this~$b$ from~$w_{i+1}$ and adding it as the penultimate letter in $y_i$ produces an RRS $V$ which results in a word with subword $ \rho_{v_i}\mathrm{pref}[\tau(\hat{v}_i)]\alpha_{v_{i+1}}\rho_{v_{i+1}}$, where this instance of $b$ ends up as $l[\rho_{v_i}]$. The commutations that allowed us to originally push $b$ to the left when we computed $\tau(v_i)$, can now be inverted to take $v$ back to $u$. } 
If, on the other hand $b$ is a pseudo-generator of $u_{i+1}$, then it does not commute with $l[\tau(\hat{u}_i)]$. In the new RRS, we have constructed~$y_i$ by inserting~$b$ into the second to last letter of~$w_i$, so $l[\tau(\hat{v}_i)]\beta_{v_i}=l[\tau(\hat{u}_i)]\beta_{u_i}b$. Thus we get $u_{i+1}=v_{i+1}$, and the remainder of the RRS is unaffected by the change. {In this case we get $u=v$.}

{To check c), we apply the algorithm in \Cref{algorithm} to $w'$. If $b$ is a generator of $u_{i+1}$, when arriving to $ba$ the algorithm applies 2.b)+1) and continues as in $U$ to obtain $V$. If $b$ is internal in $u_{i+1}$, the algorithm applies 2.b)+1) or 3.d)+1), and continues as in $U$ to obtain $V$. }

 \underline{If $a=l[w_{k}]$:} Suppose that $w_{k+1}=\epsilon$ and $b=f[\gamma]$. Because $f[\gamma]=l[\tau(u_{k})]^{-1}$, this case implies that~$a$ and~$b$ are both pseudo-generators of~$u_k$, which is impossible if $m_{ab}=2$, so this case will never occur. Assuming that $w_{k+1}\neq \epsilon$ implies that $b=f[w_{k+1}]$. We define our new RRS $V$, where $n=k$, $y_i=w_i$ for $i<k$, $y_k$ is $w_k$ with a $b$ inserted to the left of $a$ and $y_{k+1}$ is $w_{k+1}$ with its first letter removed. Note that $b$ commutes with both $a$ and $l[\tau(\hat{u}_k)]$ so $b$ commutes with both the pseudo-generators of~$u_k$. This means that when we form $y_k$ by adding an extra~$b$ we still end up with $v_k$ as a P2G critical word and this extra instance of~$b$  will end up as $l[\rho_k]$. In conducting the RRS algorithm~$V$, we may have pushed $b$ to the left via some commutations, but any of the commutations can be undone and  we can use commutations to push $b$ to the right and get~$v$ to~$u$. {To check c), we apply the algorithm in \Cref{algorithm} to $w'$. When arriving to $ba$ the algorithm applies 1) twice and continues as in $U$ to obtain $V$.}

\underline{If $ab$ is contained in $w_{k+1}$:} Recall that~$u_{k+1}$ is of the form~$c\sigma$, $c$ commutes with every letter in~$\sigma$. If~$ab$ is a prefix of~$u_{k+1}$, then $k=0$, and replacing~$ab$ with~$ba$ gives a new {optimal} RRS where $n=0$ and~$v_1$ is~$u_1$ with~$b$ removed. Otherwise, $ab$~is a subword of~$\sigma$ and every letter in~$\sigma$ commutes with the letter~$c$. Changing~$ab$ to~$ba$ does not affect this at all, and the steps of the new RRS are the same (it is also an optimal RRS). {Applying again $ba=ab$ will take $v$ to $u$.}

\underline{If $a=f[\gamma]$:} because $f[\gamma]=l[\tau(u_{k})]^{-1}$, we know that~$b$ commutes with $l[\tau(u_{k})]$.  Append~$b$ to~$w_{k+1}$, and we obtain a new RRS as desired. {In this case $\eta\neq \gamma$.}

\end{proof}

\begin{lemma}\label{Lem: tau moves do not affect RRS} Suppose that $A_S$ is a 3-free Artin group. 
Consider a word~$w$ admitting an optimal RRS $U= u_1,\dots, u_k$ with factorization $w=\mu w_1\dots w_{k+1} \gamma$. Suppose also that  $w$ has a 2-generated critical subword $\zeta$, in the generators $s$ and $t$, such that $m_{s,t}>3$. Let~$w'$ be the word obtained from~$w$ by replacing the subword $\zeta$ with the subword $\tau(\zeta)$.

\medskip\noindent
a) $w'$ also admits an RRS $V=v_1,\dots, v_{n+1}$ with factorization $w'=\nu y_1\dots y_{n+1}\eta$.

\medskip\noindent
{b) Also, if ~$\zeta$ does not intersect with~$\gamma$, then $\gamma=\eta$. In this case, if $U$ takes $w$ to $u \mathrm{suf}[\gamma]$ and $V$ takes $w'$ to $v \mathrm{suf}[\gamma]$, then we can transform $v$ into $u$ by applying only commutations and $\tau$-moves.}

\medskip\noindent
{ c) In the hypothesis of b), if $\mathrm{pref}[w]\in W$, then $V$ is also optimal. }

\end{lemma}

\begin{proof}

Once again, we will prove this by looking at all the possible cases for how $\zeta$ might intersect with the factorization $w=\mu w_1\dots w_{k+1} \gamma$.

\underline{If $l[\zeta]$ is in $\mu$:} Changes in $\mu$ do not affect the $RRS$, so replacing $\zeta$ with $\tau(\zeta)$ does not affect the RRS in $w$, {which remains optimal}. {Applying again $\tau$ to $\tau(\zeta)$ takes $v$ to $u$.}

\underline{If $l[\zeta]=f[w_1]$:} This means that $f[w_1]$ starts by a power of $l[\zeta]$, but that power has to be 1 if we want the RRS to be optimal. {To define a new optimal RRS, we apply the algorithm in \Cref{algorithm}  to $w'$. {There are two options. If $\zeta$ and $\hat{u}_1$ share the same generators $s$ and $t$, and $\alpha_{u_1}$ commutes with $\zeta$, then $y_1$ is the shortest $\{s,t\}$-P2G critical suffix of $\tau(\zeta_1)\mathrm{suf}[w_1]$ such that $l[\tau(y_1)]=l[\tau(w_1)]$ (which exists by \Cref{Lemma: non-critical P2G has RRS}) and $y_i=w_{{i-1}}$ for $i>1$. Since all the letters in $\zeta$ are pseudo-generators, when the algorithm arrives to $\zeta$, it will apply steps 1), 2.a), 3.b) and 3.a) until it will give exactly the aforementioned $y_1$.  {If $\hat{u}_1$ is not generated by $s$ and $t$, then} when the algorithm arrives at $l[\tau(\zeta)]$}, it applies 3.d) {or 3.e.iv)} and incorporates the minimal critical suffix of $\tau(\zeta)$ as $y_1$}, so $y_2=\mathrm{suf}[w_1]$, and $y_i=w_{i-1}$ for $i\geq 2$. 
{In either case, we have defined $y_1$ such that after applying $\tau$ to $y_1$ we} {obtain $v_2=w_1$} and thus the RRS can proceed from here as before. 

It remains to prove that we can transform $v$ to $u$ using commutations and $\tau$-moves. {Let $\tau(\zeta)=\overline{\zeta}\underline{\zeta}$ such that $\underline{\zeta}$ is the prefix of $y_1$ coming from~$\tau(\zeta)$. In the first case, the words $v$ and $u$ differ only in that $u$ contains $u':=\mathrm{pref}[\zeta]\mathrm{pref}[\tau(w_1)]$, while $v$ contains $v':=\overline{\zeta} \mathrm{pref}[\tau(y_1)]$ ($y_1=\mathrm{pref}[\underline{\zeta}w_1]$)}.
In the second case, the words $u'$ and $v'$ are both 2-generated geodesic words representing the same group element. By \cite[Theorem 2.4]{Holt2012}, two equivalent words geodesic words in a dihedral Artin group are related by a sequence of $\tau$-moves. Thus $u$ and $v$ are related by a sequence of $\tau$-moves and commutations as desired. {The first case works the same but using also commutations to push the letters of $\alpha_{w_1}\rho_{w_1}$ to the left in order to obtain two 2-generated geodesic subwords representing the same element. }

\underline{If $l[\zeta]$ is in {$w_1$} and $\zeta$ and  $w_1$ intersect in at least two different letters:} 
If the intersection of~$\zeta$ and $w_1$ involves both $s$ and $t$, then $s$ and $t$ are pseudo-generators of~$u_1$ or $\zeta$ is completely contained in~$w_1$ and~$s$ and~$t$ are both internal. If $s$ and $t$ are both internal, then by Lemma~\ref{lemma:beta}, $\beta_1$ is either trivial or the power of some letter, so~$\zeta$ will end up in~$\alpha_1\rho_1$, and applying~$\tau$ to~$\zeta$ does not affect the RRS. The new RRS remains optimal, because when we apply the algorithm of \Cref{algorithm} to~$w'$, $\tau(\zeta)$ will be added to~$y_1$  also as internal letters.

Suppose that $s$ and $t$ are both pseudo-generators. Let $\sigma$ be the smallest subword of $w$ containing both $\zeta$ and $w_1$, and let $\sigma'$ be the corresponding subword of $w'$ containing $\tau(\zeta)$. Applying $\tau$ to a subword of $\sigma$ does not change which commutations are possible and so both $\sigma$ and $\sigma'$ must be P2G words with $s,t$ as pseudo-generators. After we apply these commutations to $w$ we get a word with subword  $\alpha_\sigma \rho_\sigma \hat{\sigma} \beta_\sigma$. We also see that $\hat{u}_1$ must be a suffix of $\hat{\sigma}$, because $u_1$ is a suffix of $\sigma$. 


Suppose without loss of generality that $l[\tau(\hat{u}_1)]=t$. Then the word $\hat{\sigma}' t^{-1}$ is a non-geodesic word in $s$ and~$t$. Suppose that $\hat{\sigma}' t^{-1}$ is not freely reduced. In this case, we must have had $l[\zeta]=l[w_1]$, implying that $\beta_{u_1}=\epsilon$.  We can define a new RRS $V$ on $w'$ by taking $y_i=w_{i+1}$ for $1\leq i \leq k$. {We see that $V$ is also optimal by applying the algorithm in \Cref{algorithm} to $w'$ .} {To see that $u$ can be obtained from $v$ by $\tau$-moves and commutations, note that they only differ in that $u$ contains $\mathrm{pref}[w_1]$ where $v$ contains $\mathrm{pref}[y_1]$. Using commutations to push the letters of $\alpha_{w_1}\rho_{w_1}=\alpha_{y_1}\rho_{y_1}$ to the left, we obtain two equivalent 2-generated geodesic words, which are related by a sequence of $\tau$-moves thanks to \cite[Theorem~2.4]{Holt2012}.} Thus~$u$ and~$v$ are related by a sequence of $\tau$-moves and commutations as desired.

On the other hand, if~$\hat{\sigma}' t^{-1}$ is non-geodesic and  freely reduced, then by \cite[Lemma~2.3]{Holt2012}, $\hat{\sigma}'$~has a critical suffix such that applying~$\tau$ to this suffix results in a word ending in $t$.  Let~$y_1$ be the suffix of $\sigma'$ such that~$\hat{y}_1$ is the smallest suffix we obtain from \cite[Lemma~2.3]{Holt2012}. Note that $\beta_{y_1}=\beta_\sigma=\beta_{u_1}$ and so we obtain a new RRS by letting $y_i=w_i $ for $i\geq 2$. {Notice that this is the RRS given by the algorithm in \Cref{algorithm} when applied to $w'$, so it is optimal.} To see that $u$ can be obtained from $v$, note that $u$ has prefix $\mu \alpha_{u_1}\rho_{u_1} \mathrm{pref}[\tau(\hat{u}_1)] $, $v$ has prefix $\nu \alpha_{v_1}\rho_{v_1} \mathrm{pref}[\tau(\hat{v}_1)] $ and after the prefixes, the words $u$ and $v$ are identical. We can then take both of these prefixes and if necessary use commutations to move $ \alpha_{u_1}\rho_{u_1}$ to the left of any letters that were in $\sigma$ and $\alpha_{v_1}\rho_{v_1}$ to the left of any letters that were in $\sigma'$. The result is two words with identical prefixes and where the remaining suffixes are equivalent geodesic words in ${s,t}$. Again, by \cite[Theorem~2.4]{Holt2012}, two equivalent geodesic words in a dihedral Artin group are related by a sequence of $\tau$-moves.  

\underline{If $l[\zeta$] is in $w_i$ for $1<i\leq k$:} Here we divide further into subcases. Suppose first that $\zeta$ is a subword of $w_i$. If we assume that~$s$ and~$t$ are internal letters of~$u_i$, then by \Cref{lemma:beta}, $\beta_i$ is either trivial or the power of some letter, so~$\zeta$ will end up in~$\alpha_i\rho_i$, and applying~$\tau$ to~$\zeta$ does not affect the RRS. {The new RRS remains optimal, because when we apply the algorithm of \Cref{algorithm} to~$w'$, $\tau(\zeta)$~will be added to $y_i$  also as internal letters.} After some commutations, we can apply again~$\tau$ to~$\tau(\zeta)$ and take~$v$ to~$u$. Also notice that if~$s$ is an internal letter, then~$t$ has to be  also an internal letter because otherwise, by 3-freeness, $s$ would be trapped between two instances of~$t$. Finally assume that~$s$ and~$t$ are both pseudo-generators of~$u_i$. We will show that this violates the optimality of the RRS. Since $i>1$, $\zeta$ is a subword that is not a prefix of $\hat{u}_i$. By \Cref{lemma: critical_subwords}, there is a proper critical suffix of $\hat{u_i}$ containing $\zeta$, and we could have started our RRS with the first critical factor as this critical suffix. This produces an RRS with a longer head than the original RRS, violating optimality.

Now suppose that $\zeta$ is not a subword of $w_i$, so $\zeta$ intersects with both $w_i$ and $w_{i-1}$. We will show that this intersection violates our assumption that the RRS $U$ is optimal, by first showing that in this case $s$ and $t$ must be pseudo-generators of both $u_i$ and $u_{i-1}$. Assume without loss of generality that $l[w_{i-1}]=s$. Thus $s$ must be a pseudo-generator of $u_{i-1}$.

If  the intersection of $\zeta$ with $w_{i-1}$ also contains $t$, then $t$ cannot be an internal letter in $u_{i-1}$ as it does not commute with $s$, so $s$ and $t$ must be the pseudo-generators of $u_{i-1}$ and $\beta_{u_{i-1}}=\epsilon$. In this case, we see that $l[\tau(\hat{u}_{i-1})]=t^\pm$ and have $f[w_i]=s^\pm$ {(by \Cref{remark:optimality})}. This means that the word $u_i$ begins at $t^\pm s^\pm$ and  $\alpha_{u_i}=\epsilon$. Hence $s$ and $t$ must also be the pseudo-generators of $u_i$. However by optimality we must have that $\alpha_{u_i}$ contains at least one letter that does not commute with $s$ and $t$.  This gives a contradiction. 

Optimality will also be violated if~$\zeta$ and~$w_{i-1}$ intersect only in a power of a generator, say~$s^\#$. Then~$s$ is a pseudo-generator of $u_{i-1}$ and $f[\hat{u}_i]$ must be the other pseudo-generator, a letter that does not commute with~$s$, say~$r$. We also know that both~$s$ and~$t$ must appear in~$w_i$ and the intersection of~$\zeta$ with~$w_i$ must be a 2-generated word with $p+n\geq m_{s,t}-1$. This means that~$s$ is trapped in~$u_i$ and therefore must be a pseudo-generator of~$u_i$. Because $m_{s,r}>3$, $s$ must appear more than once in~$u_i$ and $l[\tau(\hat{u}_{i-1})]\beta_{u_{i-1}}$ cannot contain the letter~$s$. This in turn implies that~$w_i$ must contain a trapped copy of~$t$ as well. Thus $t=r$ and~$s$ and~$t$ must be the pseudo-generators of both~$u_i$ and~$u_{i-1}$. We now again have the case where~$s$ and~$t$ are pseudo-generators of both~$u_{i}$ and~$u_{i-1}$ and again see that $\alpha_{u_i}=\epsilon$ and optimality is violated.

\underline{If $l[\zeta]$ is in $w_{k+1}$:} If $f[\zeta]$ is also in $w_{k+1}$ then applying $\tau$ to~$\zeta$ does not affect the RRS at all {and applying again $\tau$ to $\tau(\zeta)$ takes $v$ to $u$.} {The new RRS remains optimal.}

If the intersection of $\zeta$ with $w_k$ is of the form $s^\#$, then the intersection of~$\zeta$ with $w_{k+1}$ must contain both~$s$ and~$t$. However the pseudo-generators of~$w_k$ are letters, say $s$ and $r$, where either~$r$ or~$r^{-1}$ must be~$l[\tau(u_k)]$ and~$r$ cannot commute with~$s$. This means that~$\zeta$ and~$w_{k+1}$ cannot intersect in this fashion.   

Otherwise, we know that both $s,t$ appear in some suffix of $w_k$, and so $s,t$ are the pseudo-generators of~$w_k$. This contradicts our assumption that $l[\tau(u_k)]$ commutes with any letters in~$w_k$.

\medskip

Considering all the cases we have examined up to this point, we have proven the lemma in the case where $\zeta$ does not intersect with $\gamma$. Now we will consider further cases and show that applying~$\tau$ to~$\zeta$ still produces a word that admits an RRS, although in these cases the tail of the new RRS may not be equal to $\gamma$:

\underline{If $\zeta$ is a subword but not a prefix of $\gamma$:} Changes to this part of $w$ do not affect the existing RRS at all. 

\underline{If $\zeta$ is a prefix of $\gamma$:} Suppose that $f[\zeta]=f[\gamma]=t$. In this case, we apply the $\tau$-moves from our RRS to $u_1\dots u_{k}$ to obtain a word that ends in $t^{-1}w_{k+1}$. We see that $t^{-1}w_{k+1}\zeta$ is a non-geodesic P2G generated word. This implies that $t^{-1}w_{k+1}\tau(\zeta)$ is also non-geodesic and P2G, and by Lemma~\ref{Lemma: non-critical P2G has RRS} we can continue the RRS to obtain an RRS on the word $w'$.

\underline{If $\zeta$ intersects both $\gamma$ and $w_kw_{k+1}$:} Notice that $w_{k+1}$ must be trivial in this case, because any letters in $ w_{k+1}$ both cannot be equal to $f[\gamma]$ and must commute with $f[\gamma]$, and so cannot appear in~$\zeta$. Suppose without loss of generality that $f[\gamma]=s$. Then, $l[w_{k}]= t^{\pm 1}$, so $w_k$ shares both pseudo-generators with~$\zeta$. Apply the sequence $u_1,\dots,u_{k-1}$ to~$w$. Then the subword $l[\tau(\hat{u}_{k-1})]\beta_{k-1}\cdots l[\zeta]$ (or just $f[w_1]\cdots l[\zeta]$ if $k=1$) is P2G and non-geodesic. Hence this word is still P2G and non-geodesic after applying~$\tau$ to~$\zeta$. By Lemma~\ref{Lemma: non-critical P2G has RRS} there is an RRS in this subword that continues the sequence of transformations and give us the desired RRS.


\end{proof}

\begin{lemma}\label{Lemma_Case_ab=ba} {
Suppose that $w\in W$, $m_{a,b}=2$ and, in the process of computing $\hat\phi(wab)$, we apply two RRS, one cancelling $a$ and another cancelling $b$. Suppose $wb\not\in W$. Then $\hat{\phi}(wba)$ can be obtained from  $\hat{\phi}(wab)$ via only commutations and $\tau$-moves.   }
\end{lemma}

\begin{proof}

We know that $wa\notin W$ and $wb\notin W$. We have also assumed that $\hat{\phi}(wa)b\notin W$. We will start by showing that $\hat{\phi}(wb)a\notin W$ and so computing $\hat\phi(wba)$ and $\hat\phi(wab)$ both involve two RRS, one cancelling a and another cancelling b.

Observe that the commutations and $\tau$-moves of the RRS in $wb$ only happen inside $w$, and can be used to transform $w$ into $\hat{\phi}(wb)b^{-1}$. This RRS also transforms $wa$ into $\hat{\phi}(wb)b^{-1}a$. Using Lemmas~\ref{Lem:commutation does not affect RRS} and \ref{Lem: tau moves do not affect RRS}, $\hat{\phi}(wb)b^{-1}a$ has an RRS that cancels $a$. 
 Applying one more commutation we have that $\hat{\phi}(wb)ab^{-1}$ admits an RRS cancelling~$a$, as we wanted to show.

We now show that not only do $\hat{\phi}(wba)$ and $\hat{\phi}(wab)$ have the same length, but one can be obtained from the other using commutations and $\tau$-moves. 
The key idea in this argument is that the RRS algorithm is made up of a sequence of commutations and $\tau$-moves. If we start with $wa$ and apply an RRS, we change $w$ into a new word ending in $a^{-1}$ and cancel to obtain $\hat{\phi}(wa)$. This means that the result of the commutations and $\tau$-moves in the RRS (before any cancellation) is $\hat{\phi}(wa)a^{-1} a$. The exact same commutations and $\tau$-moves, take $wb$ to ${\hat\phi}(wa)a^{-1}b$ from $wb$. Because $wb$ admits an optimal RRS that cancels $b$, we can repeatedly apply \Cref{Lem:commutation does not affect RRS} and \Cref{Lem: tau moves do not affect RRS} to get that $\hat{\phi}(wa)a^{-1}b$ also admits an optimal RRS that cancels $b$. We call this RRS $U$. Again using \Cref{Lem:commutation does not affect RRS} and \Cref{Lem: tau moves do not affect RRS}, we know that after applying $U$ and cancelling the suffix $b^{-1}b$ we get a word that can be obtained from $\hat{\phi}(wb)$ via a series of commutations and $\tau$-moves.

Because $a$ and $b$ commute, the final $a^{-1}$ in $\hat{\phi}(wa)a^{-1}b$ is not a pseudo-generator in any critical word in $U$. This means that we can apply to $\hat{\phi}(wa)b$ an optimal RRS, called $U'$, where $U$ and $U'$ are identical except for that fact the $w_{k+1}$ ends in $a^{-1}$ in $U$ and this letter is missing from $w_{k+1}$ in $U'$. By uniqueness of the optimal RRS, we know that $U'$ is the RRS applied to $\hat{\phi}(wa)b$ in order to obtain $\hat{\phi}(wab)$. It follows that the word obtained by applying the RRS U to $\hat{\phi}(wa)a^{-1}b$ and then cancelling $b^{-1}b$ is  $\hat{\phi}(wab)a^{-1}$. Putting this together with the result of the previous paragraph we get the $\hat{\phi}(wb)$ can be obtained from $\hat{\phi}(wab)a^{-1}$ via commutations and $\tau$-moves.

We now repeat this process a second time. We know that $\hat{\phi}(wb)a$ admits a unique optimal RRS that results in $\hat{\phi}(wba)a^{-1}a$. Because $\hat{\phi}(wab)a^{-1}$ can be obtained via commutations and $\tau$-moves from $\hat{\phi}(wb)$, we can repeatedly apply \Cref{Lem:commutation does not affect RRS} and \Cref{Lem: tau moves do not affect RRS} to get that {$\hat{\phi}(wb)a$} admits a unique optimal RRS that cancels the final $a$. Call this RRS $V$.  Part c) of these lemmas implies that the result of $V$ after cancellation of $a^{-1}a$, {which is $\hat{\phi}(wba)$, can be obtained from $\hat{\phi}(wab)$ via commutations and $\tau$-moves.}  
\end{proof}

\begin{lemma}\label{P3}
	$\hat\phi(w(a,b)_{m_{ab}})=\hat\phi(w(b,a)_{m_{ab}})$ for any $w\in W$, $a,b\in S$ and where $m_{ab}>2$ is the length of the relation involving $a$ and $b$.
\end{lemma}

\begin{proof}

This proof is quite technical. To illustrate all the different cases, there are examples after most of them. Consider checking the examples before checking the general argument. 
	Let $\{a,b\}=\{x,y\}=\{s,t\}$.
	We need to prove that $\hat{\phi}(w\, _m(a,b))=\hat{\phi}(w \, _m(a,b))$.
 Suppose that there is some $i$ such that $w \, _i(a,b)\not\in W$ and $w \, _{i-1}(a,b)\in W$. Let $U=u_1,\dots, u_k, u_{k+1}$ be the RRS with factorization $w \, _m(a,b) = \mu w_1\cdots w_k w_{k+1} \gamma$. 
 
\smallskip

 \noindent
\textbf{Case 1)} Suppose that $w_{k+1}$ is trivial if $k>0$, and $w_{k+1}=f[\gamma]^{-1}$ if $k=0$. 

\medskip

\noindent
\textbf{Case 1.a)} If $k=0$, the RRS consists of just a free reduction. Let $\, _{j}(a,b)$ be the biggest prefix of $\, _m(a,b)$ that can be cancelled by free reductions and let $w'$ be the word obtained from $w$ after these $j$ free cancellations. Then $\hat{\phi}(w\, _m(a,b))=\hat{\phi} (w' \, _{m-j}(x,y))$. 

 We have an optimal RRS in $w\, _m(b,a)$ that produces the cancellation of the last $j$ letters in $\,_m(b,a)$: just take the critical subword $ (b^{-1},a^{-1})_j\, {}_{m-j}(b,a)$ as the first and only factor of the RRS. If $w\, _{m-j}(b,a)\in W$, then  $\hat{\phi}(w\, _m(a,b))=\hat{\phi}(w\, _m(b,a))$. For example, for $m=5,\, j=2$,
$$\hat{\phi}(b^{-1}a^{-1}ababa)=aba, \quad \hat{\phi}((b^{-1}a^{-1}ba b)ab)=aba.$$

However, $\hat{\phi}$ is computed recursively, so we must still consider the case where $w\, _{m-j}(b,a)\not\in W$. To establish notation, suppose that the optimal RRS  on $w\, _{m-j}(b,a)$ is given by $V=v_1, \dots v_n,v_{n+1}$ with factorization $\mu \zeta_1\dots \zeta_n\zeta_{n+1}\gamma$. Notice that $v_{n+1}$ must be empty because $l[w]=a^{-1}$ cannot be a letter in $v_{n+1}$ as any letter in $v_{n+1}$ must not be equal to the cancelling letter and must commute with the cancelling letter, and therefore must be essentially different from $a$. 

\bigskip

We divide further into subcases. 

\smallskip

Suppose that $\gamma=\, _{m}(b,a)$ and $\mu \zeta_1\dots \zeta_n=w$. In this case, $\hat{v}_n$ needs to be negative because both $l[w]$ and $l[\tau(\hat{v}_n)]$ are negative. Let $\hat{v}_n$ be of the form $\,_m(s^{-1},t^{-1})\eta^-$ or $\eta^-\,(b^{-1},a^{-1})_m$. The case $\eta^-\,(b^{-1},a^{-1})_m$ with $\eta^-\neq \epsilon$ contradicts the optimality of $V$: If $n=1$, then $j=m$ and~$V$ would not be optimal because there is an RRS starting in $(b^{-1},a^{-1})_m$ cancelling all the letters in $_m(b,a)$;
 if $n>1$ and $\beta_{n-1}= \epsilon$ then we have an RRS starting in $\zeta_n$ contradicting optimality; if $n>1$ and $\beta_{n-1}\neq \epsilon$ then the subword $l[\beta_{n-1}] \cdots \zeta_n$ is P2G critical and we would have an RRS with a longer head, contradicting optimality. Now, in the case $v_n=\,_m(s^{-1},t^{-1})\eta^-$, the RRS produces the free cancellation of all letters in $_{m}(b,a)$. 

Now we must argue that $\hat{\phi}(w\, _m(a,b))=\hat{\phi}(w \, _m(a,b))$. Firstly, $j=m$ happens if and only if $\eta^-$ is trivial and it is easy to see that the result holds because the only letters that are changed are the letters that are eventually removed. For $j>m$ and $\eta^-\neq \epsilon$, we prove that after cancelling $j$ times in $w\, _m(a,b)$ we obtain a word that admits and RRS and the result of this RRS will also cancel all the remaining letters in $_{m}(a,b)$. 
Consider the RRS on $V$ that we already know applies to $w$. After applying the first $n-1$ steps we obtain $\hat{v}_n =\,_m(s^{-1},t^{-1})\eta^-$. 
Since $\eta^-$ is not trivial, $(b^{-1},a^{-1})_j$ lies in $\eta^{-}$. Then removing these $j$ letters does not affect the criticality of $v_n$, and we could apply the RRS to eliminate the remaining letters in $_{m}(a,b)$, obtaining for $\hat{\phi}(w\, _m(a,b))=\hat{\phi}(w \, _m(b,a))$.

For example, for $m=5, j=2$,
$$\hat{\phi}(b^{-1}a^{-1}b^{-1}a^{-1}b^{-1}b^{-1}a^{-1}ababa)=\hat{\phi}(b^{-1}a^{-1}b^{-1}a^{-1}b^{-1}aba)=a^{-1}b^{-1}, $$ $$ \hat{\phi}(b^{-1}a^{-1}b^{-1}a^{-1}b^{-1}b^{-1}a^{-1}babab)=a^{-1}b^{-1}.$$

\medskip

Suppose that $\gamma$ is a proper suffix of $_m(b,a)$ with more than $j$ letters (the case with $j$ letters is $w\,_{m-j}(b,a)\in W$). $V$ produces $j+j'=n(\hat{v}_n)$ free cancellations. In an argument similar to the previous case, we will show that removing $(b^{-1},a^{-1})_j$ from $w$ maintains the existence of an RRS in $w\, _{m-j}(s,t)$ that cancels $j'$ more letters in $_m(a,b)$ after the $j$ free cancellations. Notice that, since $\gamma$ has more than $j$ letters, the subword $(b^{-1},a^{-1})_j\cdots l[v_n]$ is not critical and $\hat{v}_n$ is unsigned. Suppose that~$\eta$ is trivial. Since~$j$ is maximal, the letters that complete~$v_n$, in order to make it critical, come from previous steps in the RRS. This means that $j=n(\hat{v}_n)-1$ if $\beta_{-1}=\epsilon$, and $j=n(\hat{v}_n)-2$ otherwise. Then, we can apply the first $n-1$~transformations of~$V$  to cancel $j'\in\{1,2\}$ remaining letters in $_{m-j}(s,t)$. Also, if~$\eta$ is not trivial, then $(b^{-1},a^{-1})_j$ lies in~$\eta$ and removing these letters does not affect the criticality of~$v_n$, so we can apply the RRS $V$ (without those letters) and cancel $j'$ remaining letters in $_{m-j}(s,t)$.

 For example, for $m=5, j=2$,
$$\hat{\phi}(b^{-1}a^{-1}b^{-1}b^{-1}a^{-1}ababa)=\hat{\phi}(b^{-1}a^{-1}b^{-1}aba)=aba^{-1}b^{-1}, $$ $$ \hat{\phi}((b^{-1}a^{-1}b^{-1}b^{-1}a^{-1}ba) bab)=aba^{-1}b^{-1}.$$

\medskip 

\noindent
\textbf{Case 1.b)} Now suppose that $k>0$ and $w_{k+1}=\epsilon$. If $l[w]=\{a^{-1},b^{-1}\}$ we are on the previous case. First assume that $l[w_k]=h$ essentially different from $a$ and $b$. This forces $i=1$ and $\hat{\phi}(w\,_m(a,b))=\hat{\phi}(w' \, _{m-1}(b,a))$. 
Notice that there is an RRS in $w\,_m(b,a)$ that cancels the last letter: add to the RRS $U$ (that is completely contained in $w$) the word $a^{-1}\,_{m-1}(b,a)$. We claim that, if we suppose that $w\,_{m-1}(b,a)\in W$ then this RRS has to be optimal. Suppose that there is an RRS $V=v_1, \dots v_n,v_{n+1}$ with factorization 
$w\, _{m}(b,a)=\mu \zeta_1\dots \zeta_n\zeta_{n+1} l[{}_{m}(b,a)]$. Since $a$ and $b$ do not commute, $\zeta_{n+1}$ needs to be trivial. $\zeta_n$ has at least 1 letter and has to produce $l[{}_{m}(b,a)]^{-1}$ so $a$ and $b$ are the pseudo-generators of $\zeta_n$. Also, because the relation between $a$ and $b$ is of length $m$, $\zeta_n$ contains at least $m-2$ letters. We prove now that it contains $m-1$ letters, that is, that $\,_{m-1}(b,a)$ is a suffix of $\zeta_n$: otherwise, $a$ and $b$ would be the pseudo-generators of $\zeta_{n-1}$ and $\zeta_n$, so $h$ cannot be in 
$\beta_{\zeta_{n-1}}=\epsilon$ (Lemma~\ref{lemma:beta}) and cannot be pushed to the left because by 3-freeness $a$ appears in $\zeta_{n-1}$ to the left of $h$. If $h$ lies in $\zeta_n$, then $\beta_{\zeta_{n-1}}=\epsilon$, and $v_{n-1}$ needs to produce a $b^{-1}$ to save $h$ from the trapping (no $a$ can appear to the left of $h$). But this $b^{-1}$ would cancel the first letter of $_m(b,a)$, contradicting that 
$w\,_{m-1}(b,a)\in W$. Therefore, $l[\zeta_{n-1}]= h$ and $\hat{v}_n$ has either the form 
$a^{-1}\,_{m-1}(b,a)$ or 
$a^{-1}b^\#\,_{m-1}(b,a)$. The first one is optimally produced by the RRS that we proposed at the beginning.
To study the second one, observe that, since the group is $3$-free and $\beta_{u_k}=\epsilon$, $a$ appears in $w_k$ with no instances of $b$ at the right, so $b$ cannot be a part of $\beta_{\zeta_{n-1}}$ and the case is impossible.
Thus, if $w\,_{m-1}(b,a)\in W$, then the two RRS discussed above are the RRS used to compute $\hat{\phi}(w\,_{m}(b,a))= \hat{\phi}(w\,_{m}(a,b))$. For example, for $m=4$, 
$$\hat{\phi}(w abab)= \hat{\phi}(w'a^{-1}abab) = \hat{\phi}(w'bab), $$
$$\hat{\phi}(w baba)= \hat{\phi}(w'(a^{-1}bab)a)= \hat{\phi}(w'bab).$$
 
 However, because $\hat{\phi}$ is defined recursively we, must now show that $w\,_{m-1}(b,a)\in W$. Suppose on the contrary that the first letter that can be cancelled  in $w\,_{m-1}(b,a)$ is before the final letter. Let be  $V$ an RRS with the same factorization as before cancelling this first letter. If the cancelled letter is the first letter of $_{m}(b,a)$, $\zeta_n$ would be pseudo-generated by $h$ and $b$. However, we have seen that there is no instance of $b$ to the left of the last instance of $a$ in $w_k$, what contradicts the fact that $\beta_{\zeta_n}=\epsilon$ (Lemma~\ref{lemma:beta}).
 For the other letters, we show that $l[w_k]=l[w]=h$ has to be contained in $\zeta_n$. If this is true, then $h$ could not have been pushed to the left of $\zeta_n$, contradicting that $\beta_{\zeta_n}=\epsilon$. To prove that $h$ is contained in $\zeta_n$, by 3-freeness we just need to show that $\beta_{\zeta_{n-1}}=\epsilon$ for the case in which the first cancelled letter is the second-to-last one.
 If it is not trivial, since $\zeta_n$ is pseudo-generated by $a$ and $b$, $\beta_{\zeta_{n-1}}\in\{a^{\pm 1},b^{\pm 1}\}$ by Lemma \ref{lemma:beta}. If it is $b$, then $\zeta_n$ cannot have critical form; if it is $b^{-1}$ then the letter that we are cancelling is not the first one that could be cancelled; if it is $a^{\pm 1}$, this is impossible because $h$ does not commute with $a$. Therefore, $w\,_{m-1}(b,a)\in W$. 

\smallskip

The remaining case inside Case 1.b is $l[w_k]\in \{a,b\}$. This means that $w_k$ is P2G critical with pseudo-generators $a$ and $b$. If $l[w]=h$, this letter is internal in $w_k$ and can be pushed to the left. So we can suppose that $l[w]\in\{a,b\}$. We are going to deal with the case $l[w]=a$ (the other case is symmetrical).

In $\hat{u}_k$, we now have a suffix of the form $a \,_{m-j}(a,b)$. This implies that $\hat{u}_k$ must be of the form $(s^{-1},t^{-1})_j\eta \,_{m-j}(a,b)$, and $\tau(\hat{u}_k)$ has a suffix of the form $(x^{-1},y^{-1})_j$. After applying the RRS $U$ to $w \,_{m}(a,b)$ we obtain a word where the final $j$ letters can be cancelled via free reduction. 

We will now show in $w \,_{m}(b,a)$ we can produce an RRS that cancels $j$ letters from the middle of $\,_{m}(b,a)$ and produces the same result. Let $(b,a)_l$ be the maximal suffix of $\eta $ that is an alternation. Apply the first $k-1$ transformations of $U$ to $w$, to obtain $w'$. Then there is a suffix of $w'$ of the form $(s^{-1},t^{-1})_j\eta' \,(b,a)_l$. This suffix can be completed by the first $m-j-l$ letters of $_m(b,a)$ to create a critical word. If we apply $\tau$ to this critical word, we will produce $j$ free cancellations of letters in $_m(b,a)$. To see that this cancellation produces the same result as $\hat{\phi}(w\,_m(a,b))$ we need to check that the $l$ letters at the end of $_m(b,a)$ that have not been affected by the RRS are exactly $\delta((b,a)_l)$. Just notice that $l[\delta((b,a)_l)]=l[{}_m(b,a)]$. Observe that again, this RRS is optimal, because $U$ is optimal. For example, for $m=5$, $i=3$, $$\hat{\phi}((b^{-1}a^{-1}b^{-1}aab)aba)=abb, \quad \hat{\phi}((b^{-1}a^{-1}b^{-1}ab)abab)=abb.$$

Again, because $\phi$ is defined recursively, the final step in Case 1.b, is to show that the RRS on $w \,_{m}(b,a)$ is the one used to compute $\phi$. Suppose to the contrary that there is an RRS that cancels a letter further to the left in $\,_{m}(b,a)$. This would mean that we can produce a critical word of the form $(s^{-1},t^{-1})_{j'}\eta' \,(b,a)_{l'}$ with $j'>j$. This would imply the existence of an RRS in $w\,_m(a,b)$ cancelling a letter that lies before the letter in the position $i$, having a contradiction.

\bigskip
\noindent
\textbf{Case 2)} If $w_{k+1}$ is not trivial, since $a$ and $b$ do not commute, only the first letter of $_m(a,b)$ can be cancelled by the RRS, obtaining $\hat{\phi}(w\, _{m}(a,b))= \hat{\phi}(w'\, _{m-1}(b,a))$. If $w_{k+1}$ commutes with both $a$ and $b$, by the case $m=2$, the problem reduces to study $\mu w_1\cdots w_k \, _m(a,b) $ and $\mu w_1\cdots w_k\, _m(b,a)$ and this is done in the previous case. So suppose that $w_{k+1}$ does not commute with $b$. Let $\upsilon$ be the longest prefix of $w_{k+1}$ such that $l[\upsilon]$ does not does not commute with $b$. Then, by the case $m=2$, this reduces to study the problem for $w\upsilon \, _{m}(a,b)$ and $w\upsilon \, _{m}(b,a)$. If we suppose that $w\upsilon \, _{m}(b,a)$ admits an RRS, then this RRS must satisfy the conditions of Case 1 above because $l[\upsilon]$ does not commute with $b$. 

Finally, we show that $w\upsilon \, _{m}(b,a)$ does indeed admit an RRS. Apply the first $k$~$\tau$-moves in the RRS~$U$ to the word~$w$ to obtain a word that ends in $a^{-1}\upsilon$. Now continue this RRS by considering the P2G critical word $a^{-1}\upsilon \,_{m-1}(b,a)$. This results in an RRS that cancels the final letter in $w\upsilon \, _{m}(b,a)$.

\end{proof}